\documentclass{article}
\usepackage{arxiv} 
\usepackage{natbib}
\usepackage[version=3]{mhchem}
\usepackage{balance}
\usepackage{mathptmx}
\usepackage{sectsty}
\usepackage{graphicx} 
\usepackage{array}
\usepackage{droidsans}
\usepackage{charter}
\usepackage[T1]{fontenc}
\usepackage[usenames,dvipsnames]{xcolor}
\usepackage{setspace}
\usepackage[compact]{titlesec}
\usepackage{hyperref}
\usepackage{booktabs} 
\usepackage{amsfonts}
\usepackage{bm, amssymb}
\usepackage{adjustbox}
\usepackage{bbm}

\usepackage{epstopdf}

\DeclareMathOperator*{\argmax}{arg\,max}

\usepackage[ruled,vlined]{algorithm2e}
\usepackage{amsmath}
\usepackage{authblk}

\definecolor{cream}{RGB}{222,217,201}

\title{Meta-learning for sample-efficient Bayesian optimisation of fed-batch processes}

\author[1]{Becky Langdon}%

\author[1]{%
\hspace{1mm}Gabriel D. Patr\'on}%

\author[2]{\hspace{1mm}Chrysoula D. Kappatou}

\author[2]{\hspace{1mm}Robert M. Lee}

\author[2]{\hspace{1mm}Behrang Shafei}

\author[3]{%
\hspace{1mm}Jixiang Qing}%

\author[1]{%
\hspace{1mm}Ruth Misener}%

\author[4]{%
\hspace{1mm}Mark van der Wilk}%

\author[1]{%
\hspace{1mm}Calvin Tsay\thanks{Corresponding author: \texttt{c.tsay@imperial.ac.uk} }}%

\affil[1]{Department of Computing, Imperial College London, United Kingdom}
\affil[2]{BASF SE, Ludwigshafen, Germany}
\affil[3]{School of Mathematical Sciences, Lancaster University, United Kingdom}
\affil[4]{Department of Computer Science, University of Oxford, United Kingdom}

\renewcommand{\headeright}{Langdon et al.}
\newcommand{\undertitle}{}
\renewcommand{\shorttitle}{Meta-Learning for Sample-Efficient BO of Fed-Batch Processes}

\begin{document}

\maketitle

\begin{abstract}
The optimisation of fed-batch (bio)chemical process recipes is subject to inherent, underlying, and unmeasurable fluctuations across batches, whose trajectories are difficult to model and costly to measure. Bayesian Optimisation (BayesOpt) is a powerful tool for sampling and optimisation of expensive-to-measure functions. Gaussian Processes (GPs), the surrogate models used in BayesOpt, are static, forecast poorly, and lack generalisation across experiments, limiting their applicability to time-varying batch processes with stochastic parameters, i.e., process fluctuations. This work investigates System-Aware Neural ODE Processes (SANODEP) as a meta-learning model to overcome the limitations of GPs and increase few-shot optimisation performance in BayesOpt. Using a penicillin batch production case study, we find that SANODEP outperforms GP-based BayesOpt in the low-data regime, resulting in improved objectives when few experimental runs are performed. These improvements are observed in both on- and off-distribution batches, highlighting the generalisation capabilities of SANODEP. Using this approach, batch process operators can accelerate the initial optimisation steps in BayesOpt by deploying meta-learning or optimise the process with fewer experiments when the experimental cost is high.
\end{abstract}

\section{Introduction} \label{sec:introduction} 

Increasing competitive pressures and digitalisation opportunities are inciting the use of data-driven approaches to make optimal operating decisions in the chemical process industry~\citep{daoutidis2024machine,thebelt2022maximizing,tsay2019110th}. Batch processes, where reactions are carried out over a fixed time horizon in a closed system, have received particular attention in both open- and closed-loop settings, where the challenge is to determine optimal control trajectories, or recipes, for high-value products, e.g., in the biological \citep{abdollahi2012lipid}, food \citep{patron2024economically}, energy \citep{patron2024economic}, metallurgical \citep{dering2021dynamic}, and pharmaceutical \citep{ashraf2022multiobjective} sectors. 
For biochemical manufacturing applications in particular, batch manufacturing is ubiquitous; however optimisation of batch processes remains difficult owing to inherent and expensive-to-estimate fluctuations in unmeasurable reaction rates,~\citep{BO_for_auto_lab_init_conds_Schilter} e.g., parametrised as biological growth rates or saturation constants. 
Assuming the true system corresponds to (or is well approximated by) a dynamic model with a given set of parameters, values for these parameters often vary from batch to batch~\citep{shokry2018data} owing to different suppliers, cell conditions, ambient conditions, etc. 
From the data science perspective, we refer to a batch run, corresponding to a fixed set of (unknown) values for these unmeasurable parameters, as an optimisation \textit{task}. 

Bayesian Optimisation (BayesOpt) is a powerful tool designed for global optimisation of expensive-to-measure, stochastic black-box functions and has found increasing applications in chemical engineering, including the design of experiments and the design of optimal processes and materials~\citep{humanoutofloop, tsay_paulson2024bayesianoptimizationflexibleefficient}. 
BayesOpt couples a probabilistic surrogate model with an acquisition function that quantifies the value of sampling input points, i.e., through balancing exploration and exploitation the black-box function can be optimised. 
Most commonly, BayesOpt employ a Gaussian Processes (GPs) as the surrogate model as has been applied to a broad range of problems including molecular design ~\citep{BO_for_design_of_materials_wang}, reaction and pathway optimisation, e.g. multistep synthesis \citep{Multi_Objective_BO_for_synthesis_Jorayev}, and engineering design problems such as analogue circuit design \citep{BO_for_circuit_design_Lyu}. This modelling approach is capable of handling a wide range of problem structures including continuous functions, discrete and categorical inputs, combinatorial and graph-structured search spaces, e.g. molecule graphs~\citep{xie2025bogrape}. This is of particular interest in the drive towards automated experimental design, as BayesOpt has been shown to work well in tandem with robotic hardware, e.g., initial condition problems \citep{BO_for_auto_lab_init_conds_Schilter}, and in autonomous reaction optimisation \citep{slattery2024automated}. However, GPs have limited interpretability and lack a model structure explicitly defining time-varying (i.e., dynamical) chemical processes. Moreover, GPs perform poorly when training data do not fully span the test space and struggle to generalise across tasks. In other words, standard BayesOpt treats each new optimisation \textit{task} as a new black-box function, requiring a new set of training data to train a GP, which can be cost prohibitive as information is effectively not retained from batch to batch.

Batch reactors can be modelled from first principles~\citep{fogler1999elements} as systems of ordinary differential equations (ODEs) or stochastic differential equations (SDEs)~\citep{bayes_inference_Differential_Eqs_Girolami2008-qi,bayesian_parameter_estimation}, including capturing uncertainty in parameter estimates attributed to, e.g., measurement noise or run-to-run variance. For some systems, the posterior, i.e., parameter estimate, is available in closed form, allowing for computationally efficient deployment. Unfortunately, for most systems, the posterior distribution must be estimated using computationally intensive, numerical techniques that may require many samples. 

When determining optimal batch reactor recipes, the optimisation of initial conditions to a batch process is essential.
These initial conditions often dictate the trajectories of batch states, thus influencing the process objective. Supposing that the parameters could be accurately estimated and that the structure of the dynamical system is known, an optimal control problem could be solved to determine the optimal initial conditions. For settings where the cost of collecting data is high, e.g., with high-value speciality chemicals or new pharmaceuticals, it may be impractical to perform the experimental campaigns required to build accurate dynamical models for optimisation of the initial conditions in each subsequent batch. 
Here, an ideal method would be data efficient and be able to determine the best recipe for a given task within only a few samples, or experiments. These \textit{few-shot optimisation} methods are often reliant on \textit{meta-learning}, where the model can leverage prior data/knowledge (either simulated and/or from previously observed tasks) to quickly adapt to a new task \citep{MAML, Ravi2017}. 
These approaches learn a so-called `base' model across a wide range of tasks, enabling the model generalise across tasks rather than be optimised for one particular task. 
Then, when presented with a particular \textit{target} task, the meta-learned model can perform optimisation efficiently, by effectively recalling from similar tasks seen during training. 
Meta-learning methods can exploit shared patterns of dynamical structure across tasks and allow for sample-efficient optimisation, but do not inherently address the issues arising from an intractable posterior.

With the advent of deep learning, many meta-learning models that can efficiently approximate posteriors have been proposed. 
Neural processes, or NPs, are a class of neural latent variable models that approximate the posterior, enabling optimisation in the intractable case \citep{Garnelo2018-NPs}. NPs define a probabilistic distribution over functions, analogous to a GP, but are parameterised by a Neural Network (NN) comprising an encoder-decoder architecture that encodes the distribution across tasks in latent space. Like many meta-learning models, generalisation across tasks comes at the expense of performance on individual tasks, attributed to underfitting. Various architecture adaptations have been proposed to mitigate underfitting by enriching the context, such as conditional neural processes (CNPs) and transformer neural processes (TNPs) \citep{CNPs, TNPs}. Nevertheless, these models are domain-agnostic and do not incorporate any of the structure physically inherent to the problem. Neural ODE Processes (NODEPs) extend NPs to dynamical systems, representing initial conditions and dynamic parameters, evolving the latent space to incorporate the desired ODE structure \citep{Norcliffe2021-NODEPs}. Access to the latent space allows for direct inference of the system parameters which govern the reactions, providing interpretability alongside cross-task capabilities. 

Our recent work~\citep{SANODEP} introduces system-aware neural ODE processes, or SANODEPs, which extend NODEP's context structure and training process. 
This modification introduces a \textit{trajectory-aware} meta-learning framework, enabling learning from context data comprising multiple different trajectory observations, including from arbitrary initial conditions. 
This recent work further introduces a novel acquisition function for using SANODEP in BayesOpt. In this work, we show how these modifications make SANODEP particularly suited for fed-batch process data, and we propose a comprehensive meta-learning framework based on SANODEP for fed-batch process optimisation. We specifically investigate few-shot BayesOpt on a representative fed-batch penicillin production process~\citep{bajpai1980mechanistic}, comparing SANODEP to traditional (problem-agnostic) GP implementations.  
Our computational results display the ability of SANODEP to learn over a broad distribution of parameters, effectively embedding prior knowledge from simulations or candidate models, leading to superior performance in few-shot optimisation. Overall, this work highlights the need for thoughtful design of the surrogate model and BayesOpt strategy, as well as the potential for meta-learning, to achieve efficient and robust fed-batch process optimisation. 

This work is structured as follows: Section 2 presents the background required to deploy SANODEP, Section 3 details the SANODEP algorithm and its application to fed-batch processing, Section 4 outlines the model for a penicillin production case study, Section 5 describes the results of applying  SANODEP to the penicillin plant, and Section 6 provides conclusions as to the benefit of SANODEP on fed-batch processes.

\section{Background} \label{sec:background}

For a system comprising time-dependent variables, such as concentrations, flow rates, or temperature, we define the state $\bm{x} \in \mathcal{X} \subset \mathbb{R}^{d_x}$. The derivatives of these state variables can be modelled as a system of Ordinary Differential Equations (ODEs), written as:
\begin{equation}
\begin{aligned}
    &\dot{\bm{x}} = \frac{d\bm{x}}{dt} = \bm{f}_{\bm{k}}(\bm{x}, t), \\
    &\bm{x}(t = t_0) = \bm{x}_0,
\end{aligned}
\label{eq:ode}
\end{equation}
where $t \in \tau \subset \mathbb{R}$ and $\tau= [t_0, t_\mathrm{end}]$ denote the time domain and its closed interval feasible region. The system $\bm{f}_{\bm{k}}: \mathbb{R}^{d_x} \times \mathbb{R} \rightarrow \mathbb{R}^{d_x}$ dictates the dynamics of the system, where $\bm{k} \in \mathcal{K} \subset \mathbb{R}^{d_k}$ contains the parameters that define a given task. 
In other words, an optimisation task is defined over the system \eqref{eq:ode} for some set of parameter values $\bm{k}$. 
The state at some later time $t$, $\bm{x}_t$, may be observed experimentally, or can be simulated $\bm{x}_t = \texttt{ode\_solve}(\bm{x}_0, t, \bm{f}_{\bm{k}})$, where \texttt{ode\_solve($\cdot$)} is the output of a numerical solver. These solvers are a class of numerical integration methods which compute the initial-value problem for ODEs. There are well-established explicit and implicit methods \citep{ODE_Review}: implicit methods, e.g. backward Euler, BDF, typically require solving a system of equations at each time step, often involving inversion of a (potentially large) matrix and can therefore be computationally expensive. On the other hand, explicit methods, e.g. forward Euler, classical Runge–Kutta, utilise the solution from the previous step and are cheaper per step, but may require much finer time discretisation for stability. Many implementations exist for both types of solver \citep{diffrax_kidger2021on, torchdiffeq}. Furthermore, many biochemical models are numerically unstable for common solution techniques or are stiff systems \citep{Stiff_ODEs}. Careful selection  of the ODE-solution method and tuning of step sizes are vital to reliable time integration.

We consider the task of optimising some objective function defined over the ODE system, $\bm{g}(\cdot): \mathbb{R}^{d_x}\times \mathbb{R} \rightarrow \mathbb{R} $:

\begin{equation}
    \bm{x}^*_0, t_\mathrm{end}^* = \argmax_{\bm{x}_0 \in \mathcal{X}_0, t_\mathrm{end} \in \tau}  \bm{g}(\bm{x}_{t_\mathrm{end}}, t_\mathrm{end}),
\end{equation}
or, written as a function of initial conditions:
\begin{equation}
    \bm{x}^*_0, t_\mathrm{end}^* = \argmax_{\bm{x}_0 \in \mathcal{X}_0, t \in \tau}  \bm{g}(\texttt{ode\_solve}(\bm{x}_0, t, \bm{f}_{\bm{k}})),
\end{equation}
such that the optimal initial conditions, $\bm{x}_0^*$, and optimal stopping time, $t_\mathrm{end}^*$, can be determined. Together, these decisions determine the `optimal recipe' for a batch system. The initial conditions prescribe how much of each reactant should be added and effectively dictate how the process states will evolve over time, as they represent the starting point for all batch states. Further, the stopping time defines the end-of-batch states that comprise the process product by prescribing the duration of time for which the reactions proceed. By defining states at either boundary of the time domain, we seek to optimise the objective function, which is typically a product grade or a key performance indicator (KPI) such as yield or conversion.

For the development of new processes, a mechanistic form of the ODE system \eqref{eq:ode} is often not available. We therefore model the underlying system as a black-box during the optimisation procedure. In other words, the chemical process to be optimised is a treated as a black-box function where the underlying phenomena are unknown, and the optimal batch recipe must be found using only input-output queries, i.e., experimental runs, of the system. 

\begin{figure}
    \centering
    \includegraphics[width=0.6\linewidth]{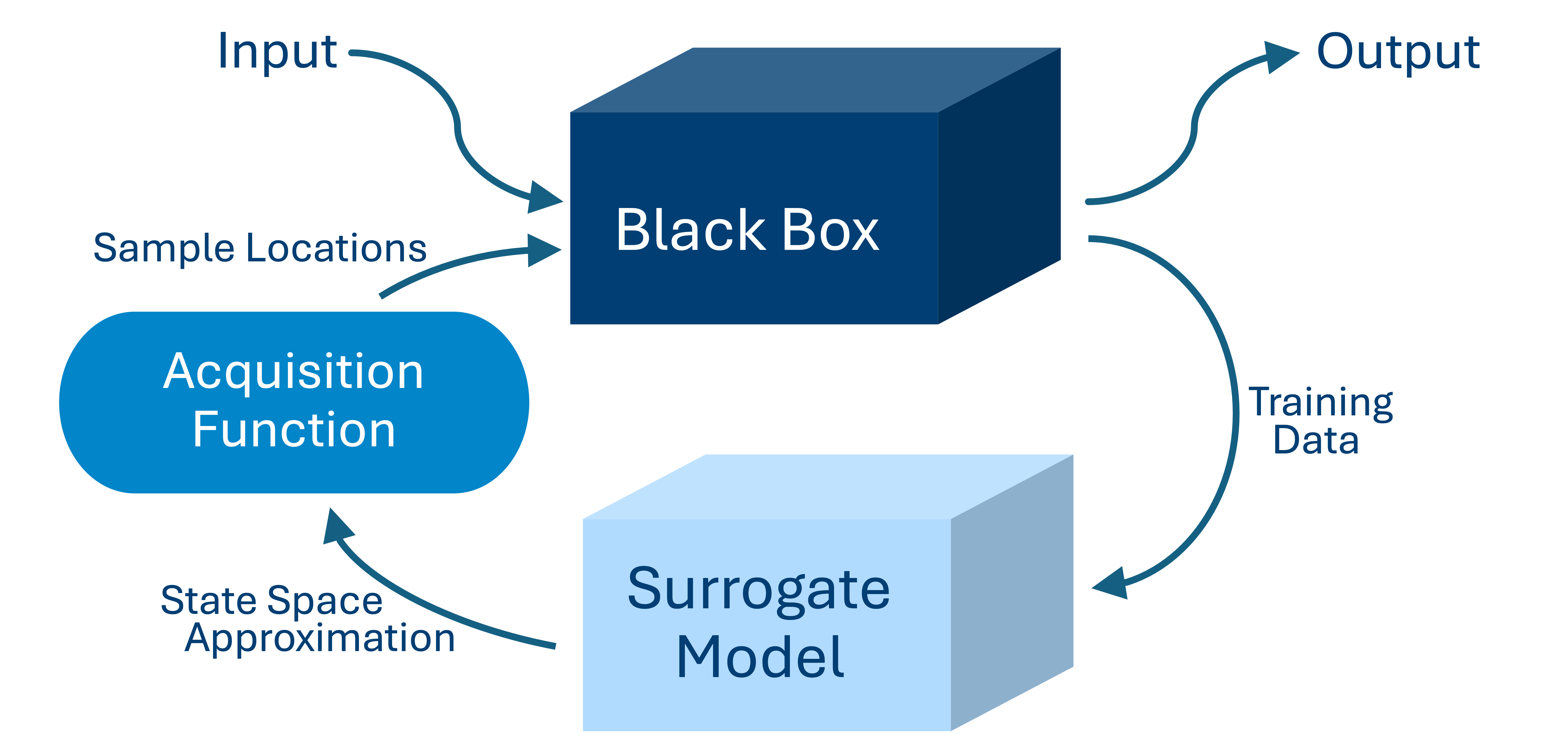}
    \caption{Conceptual depiction of the Bayesian Optimisation framework.}
    \label{fig:BO_framework}
\end{figure}

\textbf{Bayesian optimisation}, BayesOpt, provides a well established framework for black-box optimisation and is depicted in~\autoref{fig:BO_framework}. This consists of a three step process (Figure~\ref{fig:BO_framework}), which are performed iteratively until either convergence to the optimum or the experimental budget is consumed; (1) observe the black-box process, i.e., perform experiments, (2) fit a probabilistic \textit{surrogate model} to the observations, and (3) optimise an \textit{acquisition function} to determine where the black-box function should be observed next. The surrogate model must provide a robust model of the system alongside an estimation of the uncertainty for use in the acquisition function. The acquisition function then trades off exploration, minimising uncertainty in the surrogate model, and exploitation, sampling in regions which optimise the objective. 
These two features are described in the following subsections. \citet{tsay_paulson2024bayesianoptimizationflexibleefficient} review BayesOpt technologies and applications in chemical process systems.

\subsection{Probabilistic Surrogate Model} \label{sec:background_surrogate_model}

Given the requirement for a \textit{probabilistic} model, i.e., a posterior predictive distribution, traditional BayesOpt approaches usually deploy Gaussian Processes (GPs) as surrogate models \citep{frazier2018tutorial}. 

\subsubsection{Gaussian Processes}
A GP, $\mathit{GP}(\mu, \sigma): \mathcal{X} \rightarrow\mathcal{Y}$,  is a non-parametric model which defines a distribution over a family of functions, characterised by a mean, $\mu$ and a kernel (covariance) function, $\sigma$, which maps the input space $\mathcal{X}$ to the output space $\mathcal{Y}$ \citep{GP_Rasmussen_Book}. Although usually designed for non-dynamic measurements, i.e., algebraic functions, GPs have been adapted for use in dynamic processes. A common method is to treat time (or sample index) as another regressor~\citep{kocijan2013sequencing}, or as an index for multiple GP models~\citep{barton2021multivariate, zhou2015recursive}, and directly fit GP models to time profiles (also seen as signal smoothing). 

Consider a naive approach: one can adapt the GP-based BayesOpt framework to the task of black-box ODE optimisation by simply treating dynamic experiments as an input-output relationship. Specifically, one can define the input, $\bm{X}$, using the space of optimal recipe parameters $\mathbf{X} = [\bm{x}_0, t_\mathrm{end}] \in \mathcal{X}$, over which a mean, $\mu(\bm{X})$ and covariance function, $\sigma(\bm{X}, \bm{X}')$ are defined. For BayesOpt, we assume that the trajectories, or the output, $\bm{x}_t$, can be sampled from a Gaussian Process, 
\begin{align}
    \hat{f}_{\bm{k}|\bm{X}, \mathbf{\theta}_{GP}} &\sim \mathit{GP}\left( \mu(\bm{X}), \sigma(\bm{X}, \bm{X}') \right)  \\
    \bm{x}_t &= \bm{x}_0 + \int_{a=t_0}^t f_{\bm{k}}(\bm{x}, a) \, da \\
    &= \hat{f}_{\bm{k}}(\bm{X}),
\end{align}
where $\mathbf{\theta}_{GP}$ is a set of hyperparameters which parametrises the GP. In the setting presented in this paper, the function $\hat{f}$ is intended to take the form of an ODE. However, here in practice, a GP is fitted to observed data, and the resulting GP samples only approximate ODE-like behaviour. That is, the functions $\hat{f}$ produced by the GP are not strictly ODEs, even if they may resemble ODE solutions.

In defining a GP, both the prior mean and covariance functions are modelling decisions. For the mean, common choices may be a zero or constant mean, although something more expressive like a linear or quadratic mean may also be used. The covariance function then controls the deviation of samples from this mean. Common functions include the Matérn kernel, the Radial Basis Function (RBF), a special case of the Matérn kernel, and the Exp-Sine-Squared kernel~\citep{GP_Rasmussen_Book}. These functions are specified by their hyperparameters, $\mathbf{\theta}_{GP}$, referenced earlier. 
The mean and covariance are defined in relation to the GP, 
\begin{align}
    \mu(\bm{X}) &= \mathbb{E}\left[\hat{f}(\bm{X})\right], \\
    \sigma(\bm{X}, \bm{X}') &= \mathbb{E}\left[(\hat{f}(\bm{X}) - \mu(\bm{X}))(\hat{f}(\bm{X}') - \mu(\bm{X}')\right],
\end{align}
where the mean, $\mu(\bm{X})$, defines the expectation over samples from the GP, $\hat{f}(\bm{X})$, and the covariance $\sigma(\bm{X}, \bm{X}')$ the variance of these samples \citep{GP_Rasmussen_Book}. Any observed data, $\mathcal{D} = \{t_i, \bm{x}_{0, i}, \bm{x}_{t, i}\}_{i = 0}^{n_{obs}}$, can then be used to fit the hyperparameters of the mean and covariance by maximising the marginal likelihood of the GP. 

After training, i.e., once the hyperparameters have been tuned, this outlined structure allows for the posterior mean and posterior variance to be computed analytically \citep{GP_Rasmussen_Book}. This \textbf{closed-form predictive mean and variance} makes GPs well suited to BayesOpt, in comparison to more complex surrogate models, which often cannot be written in a closed form. Further, this structure allows for \textbf{explicit uncertainty quantification}. In other words, rather than returning a static prediction for each query, the GP returns a full predictive distribution for each query including a predictive mean and a predictive variance. This feature also makes GPs particularly well suited to BayesOpt, as uncertainty can be directly used in the selection of future samples via the acquisition function. We refer to this complete predictive distribution as the  \textbf{posterior}, i.e., after data are observed. 

Note that for the outlined GP, the predictive distribution does not explicitly depend on the task variable $\bm{k}$. This is because the model defines a single joint Gaussian distribution over function values indexed only by the input $\bm{X} = \left[\bm{x}_0, t\right]$. Once the GP has been fit to the data observed, the posterior collapses to a unique mean and variance at each input location. In other words, a GP cannot represent two different valid outputs for the same 
$\left(\bm{x}_0, t\right)$ pair, which would be required to \textbf{meta-learn} across multiple tasks with different underlying dynamics (different $f_{\bm{k}}$). To model task-dependent structure, an additional input dimension could be appended to input vector 
$\bm{X}$ to represent 
$\bm{k}$. This could be a learnt latent variable, as in GP-LVMs \citep{lawrence2003gaussian,titsias2010bayesian}, inferred task values, e.g., via parameter estimation, or directly measured task values. For predictions on a new and unobserved task, training data are required to determine $\bm{k}$, or a representation of $\bm{k}$, before the GP can be fit, making this approach unsuitable for few-shot optimisation. Furthermore, GPs struggle to extrapolate outside of the observed space due to their model structure. As the input moves away from observed data, the predictive uncertainty grows and the model naturally reverts to the prior distribution. Considering that during experimentation one may wish to observe the trajectory mid-run and extrapolate forwards to decide whether the objective will be optimised by continuing with the current setup or terminating early, this makes GPs further unsuitable. In short, although GPs interpolate well between observed points, their extrapolations deteriorate quickly, making predictions in these partially-observed, forward-looking scenarios unreliable. 

\subsubsection{Neural Processes}
\label{sec:NPs}
\begin{figure}
    \centering \includegraphics[width=0.35\linewidth]{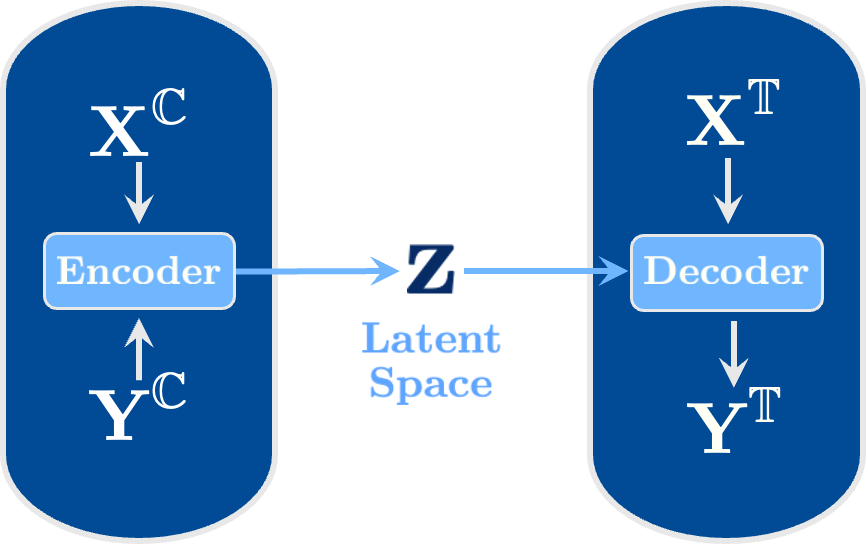}
    \hspace{10mm} \includegraphics[width=0.35\linewidth]{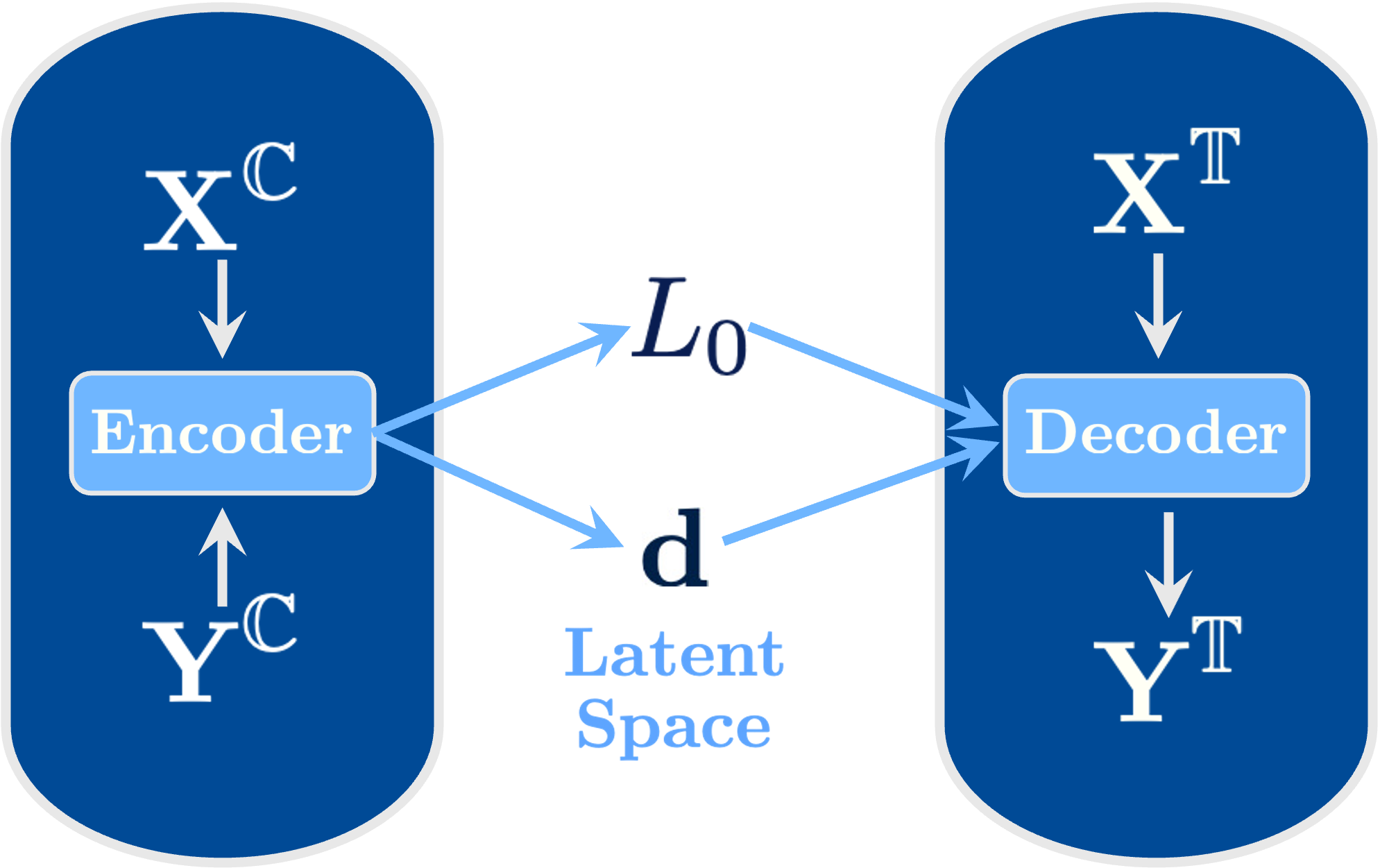}
    \caption{Left: Schematic of a Neural Process (NP). Right: Schematic of a Neural ODE Processes (NODEP). The latent space \textbf{Z} allows the NPs and NODEPs to learn across processes without an explicit representation of the task variable. }
    \label{fig:NPvsNODEP_Architecture}
\end{figure}

The challenges outlined above motivate the need for an approach that can exploit the shared structure across a family of dynamical systems, rather than treating each data batch as an isolated optimisation task, whilst still retaining the reliability and robustness of GPs. Neural Processes, NPs \citep{Garnelo2018-NPs}, are a class of models which, much like a GP, learn a distribution. However, instead of learning a distribution over a single function, an NP learns a distribution over a family of functions. This can be understood as follows: a GP can learn within a specific fed-batch process for fixed system parameters $\bm{k}$, but an NP can learn across fed-batch processes without observing the system parameters directly. An NP achieves this by encoding the observed data, or the \textbf{context dataset}, from which it infers a learnt latent representation summarising the underlying system behaviour; i.e., a representation of the task or the set of system parameters $\bm{k}$, which prescribes the fed-batch process. A schematic of this architecture can be seen in \autoref{fig:NPvsNODEP_Architecture} (left). The context set used for training may contain observations drawn from multiple fed-batch processes, and when it does, the model \textbf{meta-learns} across optimisation tasks. In other words, the model implicitly learns to make predictions across different values of system parameters $\bm{k}$. As a result of this \textbf{amortisation} over the data, the NP can make informed predictions on new and unseen tasks, i.e., tasks not present in the training set, using only a small number of observations within the new task.

Formally, this meta-learning over tasks can be expressed probabilistically; NPs \citep{Garnelo2018-NPs} model the stochastic process, $\hat{f}: \mathcal{X} \rightarrow\mathcal{Y}$, such that specific instances of the system, $\hat{f}_{\bm{k}}$, are modelled $\bm{k} \sim \pi_{k}$ where $\pi_{k}$ is a distribution over the system parameters $\bm{k}$. The context set, $\mathbb{C}$, is given $\mathbb{C} = \{\bm{X}^{\mathbb{C}}_i,  \bm{Y}^{\mathbb{C}}_i\}_{i=1}^{n_{C}}$ where $n_{C}$ is the cardinality of the training data, $\bm{X}^{\mathbb{C}}$ is the input context data given trajectory $\bm{X}^{\mathbb{C}} = \{({\bm{x}_0}, t)\}_{i=1}^{n_{C}}$, and $\bm{Y}^{\mathbb{C}}$ is the output context data $\bm{Y}^{\mathbb{C}} = \{\bm{x}_t\}_{i=1}^{n_C}$. The context captures the ``prior'' information available at training time. In our batch reaction setting, the context is the observed data from past batch processes or simulated data generated from the ODE solver, \texttt{ode\_solve}. 
The encoding distribution $p_{\theta}(\bm{z}|\mathbb{C})$ is parameterised by the encoder, a neural network (NN), $h_{e}(\mathbb{C})$, such that the latent variable is learned, $\bm{z} = h_{e}(\mathbb{C})$. At this stage, the latent variable, $\bm{z}$, is not a random variable. As in other probabilistic models, such as a Variational Auto-Encoder (VAE) \citep{Kingma_VAE}, the distribution over $\bm{z}$ can be parameterised by an intermediary, $\bm{r}$, which can be used to sample from a normal distribution inducing randomness. Here, this intermediary, $\bm{r}$, takes the form $\bm{r} = h_{e}(\mathbb{C})$. The output of this NN can then be used to sample the latent variable enforcing Gaussianity, $\bm{z} \sim \mathcal{N}\left(\bm{z}|\mu_{r}(\bm{r}), \text{diag}({\sigma}(\bm{r}))\right)$ where $\mu_r$ and $\sigma_r$ are additional encoding layers as in a VAE. 
This implicitly specifies a latent distribution $p_\theta(\bm{z}|\mathbb{C})$, where $\theta$ are the model parameters (specifically those involved in the encoder). 

The target set, $\mathbb{T}$, is constructed $\mathbb{T} = \{\bm{X}^{\mathbb{T}}_i, \bm{Y}^{\mathbb{T}}_i\}_{i = 1}^{n_{T}}$, where $n_{T}$ is the cardinality of the data for which we wish to make predictions, for the input data points, $\bm{X}^{\mathbb{T}} = \{(\bm{x}_0, t)_i\}_{i = 1}^{n_{T}}$, on which we would like to predict the outputs $\bm{Y}^{\mathbb{T}} = \{(\bm{x}_t)_i\}_{i = 1}^{n_{T}}$. 
A decoding network, $h_{d}(\textbf{X}^{\mathbb{T}}, \textbf{z})$, is trained to give predictions $\hat{f}(\bm{X}^{\mathbb{T}}) = h_{d}({\bm{X}}^{\mathbb{T}}, \bm{z})$ from the latent space. 
Because $\bm{z}$ is random, this decoder effectively defines an intermediary predictive distribution, $p_{\theta}(\textbf{Y}^{\mathbb{T}}|\textbf{X}^{\mathbb{T}}, \textbf{z})$, where $\theta$ contain the parameters of the decoder network. Note that we follow the standard NP notation, where $\theta$ effectively contains all model parameters. This distribution is subtly different from the true predictive distribution, $p_{\theta}(\textbf{Y}^{\mathbb{T}}|\textbf{X}^{\mathbb{T}}; \mathbb{C})$, we wish to sample from. However, we may  marginalise over $\textbf{z}$, 
\begin{align}
    p_{\theta}(\bm{Y}^{\mathbb{T}}|\bm{X}^{\mathbb{T}}; \mathbb{C}) = \int p_{\theta}(\bm{z}|\mathbb{C})\prod_{i = 1}^{n_{T}}p_{\theta}(\bm{Y}_i^{\mathbb{T}}|\bm{X}_i^{\mathbb{T}}, \bm{z)} \, d\bm{z}.
    \label{eq:cond_prob_NP}
\end{align}

Unfortunately, the resulting distribution is not analytically tractable but various methods exist that can be used estimate the posterior predictive distribution. In proposing the NP framework, \citet{Garnelo2018-NPs} use amortised Variational Inference (VI) to estimate this distribution. VI is a common inference technique used to bound the training loss for intractable problems through the introduction of an approximate posterior distribution, typically denoted as $q$ and often constructed $q(\bm{z}) \approx p(\bm{z}|\bm{x})$, such that (1) the distribution $q$ is much simpler to work with, and (2) the parameterisation of $q$ can be obtained by maximising a derived bound. We refer the interested reader to \citet{VI_Review} for a complete discussion of VI methods. 

The bound derived in VI is the \textbf{Evidence Lower Bound} (ELBO), which in this setting bounds the loss, $\mathcal{L} = \log p_{\theta}(\bm{Y}^{\mathbb{T}}|\bm{X}^{\mathbb{T}};\mathbb{C})$, defined as the likelihood amortised over the data, $\mathbb{D} = \{(\bm{X}^{\mathbb{C}}, \bm{Y}^{\mathbb{C}}), ( \bm{X}^{\mathbb{T}}, \bm{Y}^{\mathbb{T}})\}$, which is an aggregation of the context and target observations \citep{Garnelo2018-NPs}. This uses the amortised distributions, $q(\bm{z}|\mathbb{C})$, and $q(\bm{z}|\mathbb{D})$, as an estimator of the posterior distributions $p(\bm{z}|\mathbb{C})$, and $p(\bm{z}|\mathbb{C} \cup \mathbb{T})$, respectively. In practice, these distributions have already been learnt by the encoder, $q(\bm{z}|\cdot)$, and so the encoding network can also be used for evaluating the ELBO. The ELBO can then be derived, 
\begin{equation}
    \log p_{\theta}(\bm{Y}^{\mathbb{T}}|\bm{X}^{\mathbb{T}};\mathbb{C}) \geq \underbrace{\mathbb{E}_{\bm{z}\sim q(\bm{z}|\mathbb{D})}\left[\log p_{\theta}(\bm{Y}^{\mathbb{T}}|\bm{X}^{\mathbb{T}}, \bm{z})\right]}_{{\mathcal{L}_{E}}; \text{ Intermediary Estimation}} - \underbrace{\text{KL}(q(\bm{z}|\mathbb{D}) || q(\bm{z}|\mathbb{C}))}_{{\text{KL}_{\bm{z}}}; \text{ KL Correction}}, 
    \label{eq:NP_ELBO}
\end{equation}
which can be understood as the empirical log likelihood of the intermediary predictive distribution, $\mathcal{L}_{E}$, corrected by the Kullback-Leibler (KL) divergence term, $\text{KL}_{\bm{z}}$, which evaluates the difference between the distributions: the amortised VI inference network $q(\bm{z}|\mathbb{D})$ and the encoding of the context, $q(\bm{z}|\mathbb{C})$. Samples from the amortised posterior, $\bm{z}\sim q(\bm{z}|\mathbb{D})$, can be used to estimate $\mathcal{L}_{E}$ using typical particle methods, e.g. Monte Carlo sampling. The $\text{KL}$ divergence can be calculated analytically if, as discussed earlier, $\bm{z}$ is parametrised as a Gaussian. 
Note that an alternative to the VI framework is to simply estimate the loss, $\log \text{p}_{\theta}(\bm{Y}^{\mathbb{T}}|\bm{Y}^{\mathbb{T}}; \mathbb{C})$ in \autoref{eq:cond_prob_NP}, directly using Monte Carlo methods \citep{ConvNP}. In either case, to then train the NP, these losses can be optimised using standard gradient-based optimisation~\citep{wright1999numerical}. 

Various architectures have been proposed for the encoding and decoding networks including Convolutional Neural Processes (CNPs) \citep{ConvNP}, which incorporate a Convolutional Neural Network (CNN) in both the encoder and decoder at the cost of consistency in $\bm{z}$, Attentive Neural Processes (ANPs) \citep{Kim_ANP}, which incorporate a self-attention mechanism to address the under-fitting problem present in NPs, and Transformer Neural Processes (TNPs) \citep{TNPs}, which extend ANPs to sequential decision modelling problems using Transformer architecture. The choice of architectures is a key modelling decision which enables the surrogate model to better adapt to the problem space. See \citet{VI_Review} for a modern review of VI and \citet{NP_Family_GitHubRepo} for an accessible implementation of the models and inference methods discussed. 

\subsubsection{Neural ODE Processes}
Of particular interest to the ODE setting presented herein are Neural ODE Processes (NODEPs)  \citep{Norcliffe2021-NODEPs}, which incorporate ODE structure directly into the NP model through both a decomposition of the latent space into meaningful representations and explicit ODE structure in the latent space, see \autoref{fig:NPvsNODEP_Architecture} (right). Within the encoder, the latent structure is adapted to the problem by decomposing the latent variable, $\bm{z}$, into two separate latent variables: the latent initial conditions, $\bm{l}(t_0) \in \mathbb{R}^{n_l}$, and the time-independent, latent control signal, $\bm{d} \in \mathbb{R}^{n_{\bm{\bm{d}}}}$. This decomposition allows for meaningful, physical interpretation of the latent space, but it is worth noting that the latent variables do not have to explicitly model the task parameters $\bm{k}$ or initial conditions $\bm{x}_0$. 

At this stage, a modelling decision must be made regarding how to parametrise these latent variables. As discussed above, a simple assumption is to treat $\bm{l}(t_0)$ and $\bm{d}$ as Gaussian random variables, parametrised together by a vector $\bm{r}$~\citep{Norcliffe2021-NODEPs}. 
The parameters $\bm{r}$ are modelled using a neural network, $\bm{r} = \phi_{av}\left(\left\{h_e((\bm{X}^{\mathbb{C}}_i, \bm{Y}_i^{\mathbb{C}}))\right\}_{i=1}^{n_{\mathbb{C}}}\right)$ where the context is given for the input, $\bm{X}^{\mathbb{C}} = \{t^{\mathbb{C}}\}_{i=1}^{n_{\mathbb{C}}}$, and the output, $\bm{Y}^{\mathbb{C}} = \{\bm{x}_t^{\mathbb{C}}\}_{i=1}^{n_{\mathbb{C}}}$, for some averaging function, $\phi_{av}(\cdot): \mathbb{R}^{n_r \times n_{\mathbb{C}}} \rightarrow \mathbb{R}^{n_r}$. Note that, unlike the GP and NP models presented, the inputs $\bm{X}^{\mathbb{C}}$ and $\bm{X}^{\mathbb{T}}$ for NODEP are constructed from the time observation solely. The averaging function, $\phi_{av}$, is introduced to ensure consistent initial conditions of the latent system variable $\bm{l}(t_0)$ for context samples taken from the same trajectory but drawn at different times.  
From the global representation, the encoding network $p_{\theta}(\bm{d}|\mathbb{C})$ for the system latent variable $\bm{d}$ can be constructed $p_{\theta}(\bm{d}|\mathbb{C}) \approx q_{D}(\bm{d}|\mathbb{C}) = \mathcal{N}(\bm{d}|{\mu}_D(\bm{r}), {\sigma}_D(\bm{r}))$, where $\mu_D$ and $\sigma_D$ are both additional layers in the encoder, as in a VAE and an NP. The second encoding network, $p_{\theta}(\bm{l}(t_0)|\mathbb{C})$ for the latent initial conditions variable $\bm{l}(t_0)$ uses the initial conditions $\bm{x}_0$ as they are directly available such that $p_{\theta}(\bm{l}(t_0)|\mathbb{C}) \approx q_{L}(\bm{l}(t_0)|\bm{x}_0) = \mathcal{N}(\bm{l}(t_0)|{\mu}_L(\bm{x}_0), {\sigma}_L(\bm{x}_0))$, where again $\mu_L$ and $\sigma_L$ are both additional layers as outlined above.

Similarly to the GP, the conditional samples generated from an NP are only approximations of ODE functions. On the other hand, NODEPs enforce explicit ODE structure, as in a Neural ODE (NODE) \citep{chen2019_NODE}, by evolving the latent space through time. The latent ODE space is evaluated on a schedule of times given in the target input, $\bm{X}^{\mathbb{T}} = \{t_i\}_{i = 0}^{n_{\mathbb{T}}}$. The latent state space is evolved pointwise, 
\begin{equation}
    \bm{l}(\bm{X}^{\mathbb{T}}_i) = \bm{l}(t_0) +  \int_{a=t_0}^{t^{\mathbb{T}}_i} h_{\theta_{ODE}}(\bm{l}(a), \bm{d}, a) \, da,
\end{equation}
where $h_{\theta_{ODE}}$ is a neural network that parameterises the latent dynamics. In practice, an ODE solver, \texttt{ode\_solver}, can be used to solve the latent ODE problem, though recent methods have explored full discretisation of the time domain for NODEs~\citep{SHAPOVALOVA2025469}. Finally, as in an NP, a decoding network, $h_d(\bm{l}(\bm{X}^{\mathbb{T}}), \bm{X}^{\mathbb{T}})$ is introduced to parametrise the posterior predictive distribution $p(\bm{Y}^{\mathbb{T}}|\bm{l}(\bm{X}^{\mathbb{T}}), \bm{X}^{\mathbb{T}})$ as a function of the latent space. \citet{Norcliffe2021-NODEPs} present the full training loss, derived using amortised Variational Inference as in an NP.  These adaptations to the encoder architecture allow NODEPs to make interpretable and physically consistent predictions on unseen systems using limited task data. 

The choice of surrogate model is shaped by the structure of the domain and the optimisation problem, and no single model is universally optimal. GPs offer strong theoretical guarantees, but their ability to generalise beyond the observed region is limited, making them less suitable when data are sparse or when extrapolation is required. In contrast, NPs can meta-learn across related tasks, enabling rapid adaptation and few-shot learning in low-data regimes. However, this flexibility comes with its own challenges, e.g., underfitting.  
Once a surrogate model has been selected, the next step in establishing the Bayesian optimisation framework is the choice of acquisition function. This component leverages both the predictive mean and predictive uncertainty of the chosen surrogate to guide sampling and ensure efficient experimentation.

\subsection{Acquisition Function} \label{sec:background_acquisition_function}

Within the BayesOpt framework, the acquisition function determines where to observe the black box next. This is done by finding which set of inputs to the surrogate model optimises the acquisition function, which can be a challenging optimisation problem in itself~\citep{xie2025global}. The output of this stage is a location, or series of locations if optimising a schedule, at which the black-box process should be queried. In practice, this sampling stage corresponds to a `design of experiments,' taking measurements from a physical experiment or sampling from an expensive computational approximation. This additional ground truth data would then be passed back into the surrogate model for retraining. It is therefore vital, that the acquisition function both explores the space, to avoid local optima, and maximises the optimisation target. It is important to note at this stage that, the output of the surrogate models must be \textbf{uncertain}, i.e., Bayesian, and provide a distribution over the predictions such that exploration can occur. 

The selection of an acquisition function is a vital modelling decision which impacts the speed and reliability of convergence of on the optimal solution \citep{Wilson_BO_Acq}. There exists a broad range of acquisition functions; some of the most commonly used include Upper Confidence Bound (UCB), Expected Improvement (EI), Probability of Improvement (PI) ~\citep{frazier2018tutorial}.

Suppose a candidate surrogate model has been selected. This surrogate could either model the trajectories, $\tilde{f}(\cdot): \mathbb{R}^{d_x} \times \mathbb{R} \rightarrow \mathbb{R}^{d_x}$, to be then used in computing the optimisation objective function $g(\cdot)$, or model the objective directly, $\tilde{g}(\cdot): \mathbb{R}^{d_x} \times \mathbb{R} \rightarrow \mathbb{R}$. For the data available, $\mathbb{D} = \{(\bm{x}_0, t, \bm{x}_t)_i\}_{i=1}^{n_{\mathbb{D}}}$, the optimal observed value can be determined. Depending on the optimisation problem, this would be either the maximum, $g_{\max} = \max\left(g(\bm{x}_t)^{(1)}, \dots, g(\bm{x}_t)^{(n_{\mathbb{D}})}\right)$, or the minimum, $g_{\min} = \min\left(g(\bm{x}_t)^{(1)}, \dots, g(\bm{x}_t)^{(n_{\mathbb{D}})}\right)$. For a problem where we wish to maximise the objective function, the improvement is defined $I(\bm{x}_0, t) = \max\left( \tilde{g}(\bm{x}_0, t) - g_{\max}, 0 \right)$. If the predicted value of the surrogate model is likely to be less than the best observed value, the predicted improvement is therefore low. Following this notion, the acquisition function, \textbf{Expected Improvement}, $\text{EI}$,  is defined,

\begin{equation}
     \text{EI}(\bm{x}_0,t) = \mathbb{E}_{\bm{x}_0\in \mathcal{X},t\in\tau}\!\left[
     \max\!\bigl(0, g(\bm{x}_0,t)-g_\mathrm{end}\bigr)\right],
\end{equation}
as the expectation over the input state space \citep{Jones_EI_98}. 

EI in combination with a GP surrogate model is widely used across a broad range of chemical applications \citep{gp_bo_example_chromatography, gp_bo_example_catalysts, gp_bo_example_continous_flow}. For an NP,  a similar expression can be found and evaluated for EI and used as an acquisition function in BO \citep{NP_BO}.

 \section{System Aware Neural ODE Processes} \label{sec:model}

BayesOpt is a powerful optimisation tool, and its usefulness is further enhanced when the surrogate model reflects domain-specific structure. Meta-learning offers a way to exploit the shared structure that exists across a family of related optimisation tasks, rather than treating each batch as an entirely new problem. In many experimental BayesOpt settings, batches differ in their kinetics or ambient conditions but still follow broadly similar dynamical patterns. A model that can internalise these recurring patterns can adapt far more quickly when faced with a new batch, reducing the data required to make reliable forecasts and improving performance in the early stages of optimisation.

With this motivation in mind, we present System Aware Neural ODE Processes (SANODEP): a novel surrogate modelling approach capable of meta-learning dynamic trajectories. The SANODEP framework is described in our preliminary proceedings\citep{SANODEP} and is presented below. SANODEP builds on NODEP by incorporating an adaptable loss function and training process to allow the surrogate model to effectively switch between two subtly different optimisation modes: forecasting and interpolation. Forecasting occurs when a new run of the experimental system begins and the only data point available for this trajectory is the initial conditions. Interpolation, on the other hand, occurs when one or more measurements have been taken during the run. In both cases we wish to make model predictions on the system state at a schedule of later times. This flexibility enables SANODEP to meta-learn across heterogeneous systems, improving the ability to generalise across tasks in low-data regimes. Key to this process is the selection/generation of a meta-training dataset, which effectively captures prior information about the system. 

\subsection{Model Structure} \label{sec:model_model_structure}

Building on NODEP, SANODEP incorporates both the observation time $t$ and the initial conditions $\bm{x}_0$ into the input, $\bm{X} = \{t_i, \bm{x}_{0,i}\}_{i=1}^{n_{\text{obs}}}$, where the output makes predictions on the state at the observation time $\bm{Y} = \{\bm{x}_{t,i}\}_{i=1}^{n_{\text{obs}}}$. We refer to observations with the same initial conditions as a \textbf{trajectory}, $\mathcal{T}$, given for the $i^{th}$ trajectory of size $n_{T_i}$,
\begin{align}
    \mathcal{T}_i &= \{t_j, \bm{x}_{0,j}, \bm{x}_{t,j}\}_{j=1}^{n_{T_i}}, \\
                &= \{\bm{X}_i, \bm{Y}_i\}.
\end{align}

The context dataset, $\mathbb{C}$, is then constructed from $M$ of these trajectories $\mathbb{C} = \{\mathcal{T}^{\mathbb{C}}_i\}^M_{i=1}$. For some target dataset on which we wish to make predictions $\mathbb{T}$ we similarly construct $\mathbb{T} = \{\mathcal{T}^{\mathbb{T}}_i\}$ using the inputs to predict on $\bm{X}^{\mathbb{T}}$ and the corresponding predicted outputs $\bm{Y}^{\mathbb{T}}$. At this stage, we place no restrictions on the target. Further discussion of the design of the target in context of the training procedure can be found in \autoref{sec:model_loss_function}.

As a type of NP, an encoder is required to `encode' the physical/measured variables into the latent space, and a corresponding decoder to recover the measured variables from the latent variables. In essence, this approach seeks to represent the distribution over tasks over a set of learned latent variables, i.e., a nonphysical set of underlying variables that define the optimisation tasks. See \autoref{sec:background_surrogate_model} for further details.

SANODEP adopts the decomposition of the latent variable $\bm{z}$ into the latent initial conditions $\bm{l}(t_0, \bm{x}_{0,i}) \in \mathbb{R}^{n_l}$ for the $i^{\text{th}}$ trajectory, and the latent control signal $\bm{d} \in \mathbb{R}^{n_d}$ as in an NODEP. The encoder is used to aggregate over the trajectories in the context to produce a pair of amortised variational posteriors over the latent variables representing both the initial state and system dynamics $q_L(\bm{l}(t_0, \bm{x}_{0, i})\mid \mathbb{C})$, and $q_D(\bm{d} \mid \mathbb{C})$ respectively. These latent variables are used to parametrise the latent trajectory $\bm{l}(t, \bm{x}_{t, i})$, which is evolved as an ODE,
\begin{equation}
\bm{l}(\bm{X}^\mathbb{T}_i) = \bm{l}(\bm{X}_{0, i}^{\mathbb{T}}) + \int_{a=t_0}^{t^\mathbb{T}_i} h_{\theta_{ODE}}(\bm{l}(a, \textbf{x}_{0, i}), \bm{d}, a,  \bm{x}_{0, i}) \, da,  
\label{eq:ODEevolve}
\end{equation}
where $h_{\theta_{ODE}}$ is a neural network which parameterises the latent trajectory and $\bm{X}_{0, i}^{\mathbb{T}}$ is the target input space. Once solved, the latent state is decoded into the observed state space using a decoder network $h_{\theta_{d}}$,

\begin{equation}
\bm{Y}_i^{\mathbb{T}} \sim \mathcal{N}\big(h_{\theta_d}(\bm{l}(\bm{X}_i^{\mathbb{T}}), \bm{X}_i^{\mathbb{T}}), \sigma^2 \bm{I}\big).
\end{equation}

The full predictive distribution is therefore factorised as:
\begin{equation}
p(\bm{Y}^{\mathbb{T}}, \bm{l}_0(t_0), \bm{d} \mid \bm{X}^{\mathbb{T}}, \mathbb{C}) = p(\bm{l}_0(t_0) \mid \mathbb{C})\, p(\bm{d} \mid \mathbb{C})\prod_{i=1}^{n_T} p\big(\bm{Y}^{\mathbb{T}} \mid h_{\theta_d}(\bm{l}(\bm{X}_i^{\mathbb{T}}), \bm{X}_i^{\mathbb{T}})\big). 
\end{equation}

\begin{figure}
    \centering
    \includegraphics[width=0.4\linewidth]{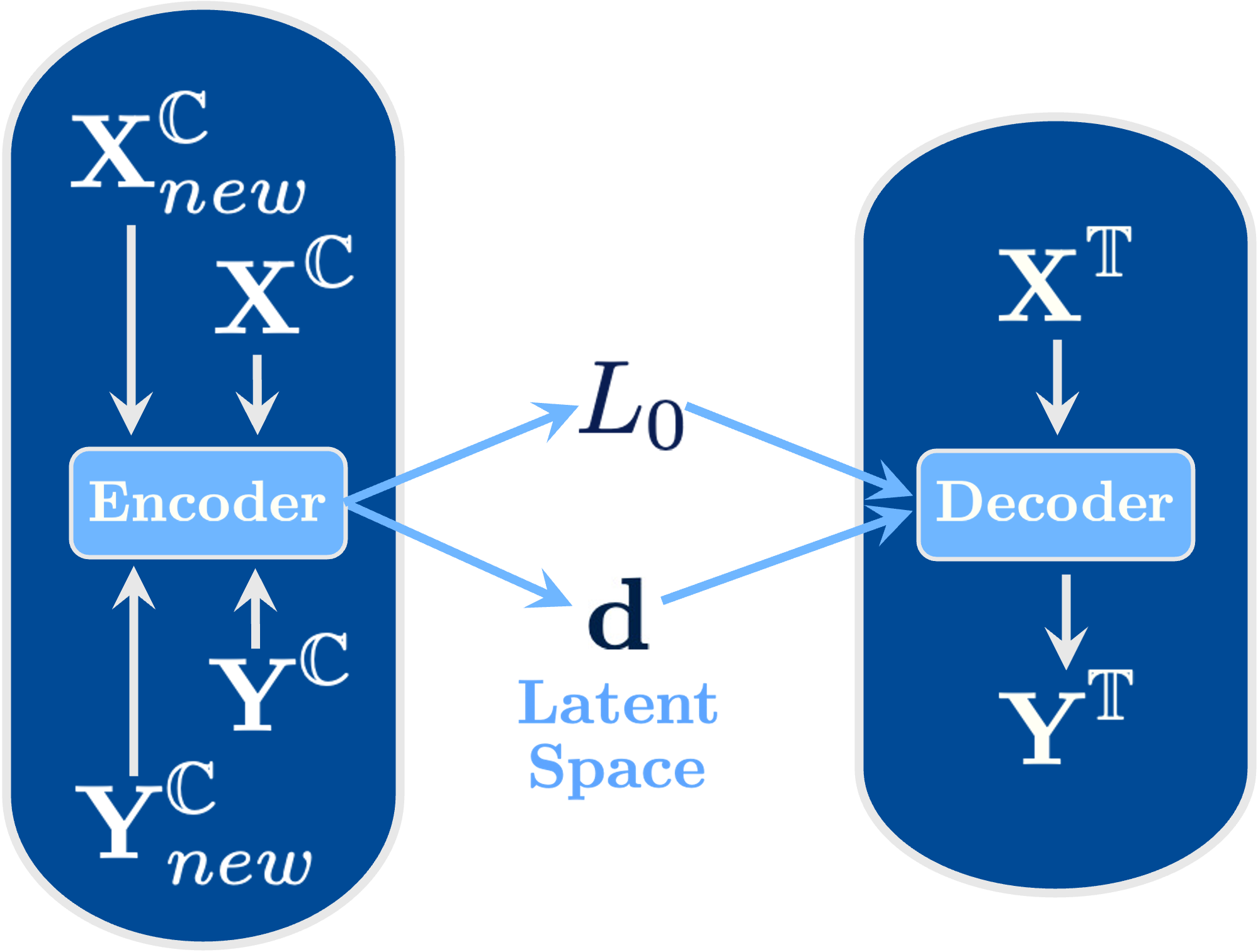}
    \caption{Schematic of the SANODEP architecture. The model receives a context set comprising previous trajectories plus new data points (forecasting or interpolation) and outputs a predictive distribution for the target set.}
    \label{fig:sanodep_architecture}
\end{figure}

Given target input $\bm{X}^{\mathbb{T}}$, we can sample the predictive output $\bm{Y}^{\mathbb{T}}$ from the distribution once the neural networks have been trained. In summary, SANODEP defines a probabilistic, surrogate model that allows for meta-learning and explicit handling of system dynamics. See  \autoref{fig:sanodep_architecture} for an overview. We next describe how SANODEP is trained for use in forecasting and interpolation modes, and can be applied to the BayesOpt framework. 

\subsection{Bi-Scenario Loss Function} \label{sec:model_loss_function}

The bespoke training procedure for SANODEP is explicitly tailored to the BayesOpt setting. The model employs a \textbf{bi-scenario} loss function designed to promote accurate predictive performance during the BayesOpt loop. This is achieved by augmenting an initial context set with new trajectories through Bayesian updates. The model outlined above is trained using episodic learning~\citep{vinyals2016matching}. Specifically, the training process is structured through multiple episodes, where each episode is designed to ensure good model performance under one of the two predictive modes; forecasting, and interpolation. A Bernoulli indicator $\mathbbm{1}$ is sampled at each episode to randomly determine the scenario on which to train (forecasting vs interpolation):

\begin{equation}
    \mathbb{C}_{\text{update}} =  \begin{cases}
         \mathbb{C}_{\text{forecast}} \; \text{ if } \mathbbm{1} = 1, \\
         \mathbb{C}_{\text{interp}}  \;\;\; \text{ else},
    \end{cases}
\end{equation}
such that a bi-scenario loss function can be implemented to train SANODEP. 
We assume that the meta-learning dataset is created using a (set of) candidate ODE model(s); effectively, the task of stating the prior distribution becomes simply expressing candidate models and parametrisations. With this assumption, the trajectories must be created arbitrarily (compared to using a historical dataset where they would be given). 
At the start of each episode an instance of the dynamical system, $f_{\bm{k}}$ is randomly sampled $\textbf{k} \sim \mathcal{K}$ alongside the initial number of observed trajectories $M$, the number of observed context points, the initial conditions of each trajectory and the observed times. An numerical ODE solver can then be used to construct the initial meta-training dataset $\mathbb{C}$ dependent on the scenario.

At the start of each episode, an initial set of $M$ observed trajectories is randomly sampled from the meta-training dataset and used to construct the observed context set, $\mathbb{C}_{\text{observed}} = \{\mathcal{T}^{\mathbb{C}}_i\}_{i=1}^M$. During the training episode, this initial context set is then augmented with further observations that the model must learn from. This information update is dependent on the predictive mode:

\begin{enumerate}
    \item \textbf{Forecasting:} In this case, the update corresponds to the selection of the ${M+1}^{\text{th}}$ trajectory, i.e., a new set of initial conditions with predictions being made at some later times on the ${M+1}^{\text{th}}$ trajectory. The context update is given,
    \begin{equation}
        \mathbb{C}_{\text{forecast}} = \{t_{0}^{M+1}, \bm{x}_0^{M+1}, \bm{x}_0^{M+1}\}.
    \end{equation}

    \item \textbf{Interpolation:} The update consists of an additional $n_o>1$ observations are made on the $i^{\text{th}}$ trajectory from the existing set of observed trajectories. The target is then further observations of this updated trajectory. 
    \begin{equation}
    \mathbb{C}_{\text{interp}} = \{t_j^{i},  \bm{x}_0^{i}, \bm{x}_j^{i}\}_{j = 1}^{n_{\text{o}}}.
    \end{equation}
\end{enumerate}

The full context set is then $\mathbb{C} = \mathbb{C}_{\text{observed}} \cup \mathbb{C}_{\text{update}}$Í where the Bernoulli indicator is used to select the updated context from the two modes. 
In effect, this training procedure forces the model to become good at learning both from new trajectories, and from new samples within existing ones. 
The full target set is then constructed in a supervised manner from existing, known observations and the desired predictions, i.e., the full context and the targets from the new trajectory: $\mathbb{T} = \mathbb{C}_{} \cup \mathbb{T}_{update}$. Unlike NODEP, which predicts on a single trajectory, SANODEP leverages information from previous tasks in each new trajectory prediction, allowing amortisation across tasks.

The flexible loss function is then given,
\begin{equation}
\mathcal{L}_\theta = 
\mathbb{E}
\left[ \log p_\theta \big(\bm{Y}_{\text{update}}^{\mathbb{T}} \mid \mathbb{C} \cup \mathbb{C}_{\text{update}}(\mathbbm{1}), \bm{X}_{\text{update}}^{\mathbb{T}}\big) \right],
\label{eq:Liklihood_SANODEP}
\end{equation}
where $p_\theta(\cdot)$ represents the SANODEP predictive distribution parameterised by $\theta$. We refer the reader to \citet{SANODEP} for further details on the computation of $\theta$. Experimentally, directly learning the task parameters, or a representation, for every new trajectory would require solving a separate, expensive problem each time. Instead, as in an NPs, discussed in \autoref{sec:NPs}, SANODEP amortises this cost by training an encoder network that learns an approximate posterior distribution over the latent variables directly from the context data, so that inference for each new task is fast, consistent, and does not require re-optimising from scratch, e.g., solving an estimation from the latent signal $\bm{d}$). The result is that the log-likelihood in Eq. \ref{eq:Liklihood_SANODEP} can be approximated via an evidence lower bound (ELBO), analogous to \eqref{eq:NP_ELBO}:
\begin{multline}
\log p_\theta(\bm{X}_{\text{new}}^{\mathbb{T}} \mid \mathbb{C} \cup \mathbb{C}_{\text{update}}, \mathbb{T}_{\text{new}}) 
\approx
\mathbb{E}_{q(\bm{d} \mid \mathbb{C} \cup \mathbb{C}_{\text{update}} \cup \mathbb{T}_{\text{new}}) \, q(\bm{l}_0^{\text{new}} \mid t_0^{\text{new}}, \bm{x}_0^{\text{new}})} 
\big[ \log p_\theta(\bm{X}_{\text{new}}^{\mathbb{T}} \mid \mathbb{T}_{\text{new}}, \bm{d}, l_0^{\text{new}}) \big] \nonumber \\
\quad - \mathrm{KL}\big( q(\bm{d} \mid \mathbb{C} \cup \mathbb{C}_{\text{update}} \cup \mathbb{T}_{\text{new}}) \,\|\, q(\bm{d} \mid \mathbb{C} \cup \mathbb{C}_{\text{update}}) \big)  - \mathrm{KL}\big( q(\bm{l}_0^{\text{new}} \mid t_0^{\text{new}}, \bm{x}_0^{\text{new}}) \,\|\, p(l_0^{\text{new}}) \big).
\end{multline}

The ELBO can then be used as the training objective for SANODEP in the bi-scenario framework. The training process is described in Algorithm~\ref{alg:sanodep}, adapted from \citet{SANODEP}.

\subsection{Bayesian Optimisation} \label{sec:model_bayes_op}

The SANODEP framework is presented for use as a surrogate model within a few-shot Bayesian optimisation framework, allowing for the optimisation of initial conditions and stopping times under limited evaluations on the target task. By leveraging the meta-learned latent dynamics, SANODEP generalises well in the low-data scenario. In other words, we leverage prior/historical knowledge (encoded in the meta-learning training) to reduce the number of experiments required to optimise a new target system.  In brief, the non-myopic strategy can be understood as:

\begin{itemize}
    \item Step 1: optimise both initial conditions $\textbf{x}_0^*$ and the measurement schedule $\textbf{t}^* = \{t_0, t_1, ...t_N\}$ to maximise an acquisition function $\alpha$ based on the SANODEP model.
    \item Step 2: acquire new observations at $t_n$ using either real-world inputs or simulated and update the context set. Then re-optimise in the remaining time window for the batch $[t_n + \Delta t, t_{stop}]$. 
\end{itemize}

In this loop, we implement a a minimum time delay $\Delta t$ between samples to enable \textbf{data-efficient optimisation}, while sitll reflecting real-world experimental constraints. 
The acquisition function is based off of Batch Hypervolume Improvement which extends Expected Improvement to the multiobjective, higher dimensional case \citep{daulton2020differentiable}. We refer the interested reader to~\citet{SANODEP} for the further details on the acquisition function formulation.

\begin{algorithm*}
\caption{Training Process in System-Aware Neural ODE Processes (SANODEP) \citep{SANODEP}}
\label{alg:sanodep}
\DontPrintSemicolon

\KwIn{Distribution of ODE systems $\mathcal{F}$, trajectory range $[M_{\min}, M_{\max}]$, 
batch sizes $N_{x_0}$ and $N_{\text{sys}}$, 
time grid $T_{\text{grid}} = \text{linspace}(t_0, t_{\max}, N_{\text{grid}})$, 
context limits $(m_{\min}, m_{\max})$, 
target limits $(n_{\min}, n_{\max})$, 
initial condition space $\mathcal{X}_0$.}

\KwOut{Trained SANODEP parameters $\theta$.}

\BlankLine
\textbf{Initialise:} Model parameters $\theta$ with random seeds.\;

\For{each training step}{
    \tcp{Generate training data}
    \For{$j = 1$ \KwTo $N_{\text{sys}}$}{
        Sample \textbf{ODE system} $f \sim \mathcal{F}$ and $N_{x_0}$ initial conditions from $\mathcal{X}_0$.\;
        \textbf{Solve ODEs} on $T_{\text{grid}}$ to calculate trajectories, $\{\mathcal{T}_i\}_{i=1}^{N_{x_0}}$.\;
        Sample number of observed trajectories, $M \sim \text{Uniform}(M_{\min}, M_{\max})$.\;

        \For{$l =1$ \KwTo $M$}{
        Randomly \textbf{subsample Target set}, $\mathcal{T}_l^{T} \subseteq \mathcal{T}_l$, s.t. $|\mathcal{T}_l^T| = n_l$ where $n_l \sim \text{Uniform}(n_{min}, n_{max})$.
        
        Randomly \textbf{subsample Context set}, $\mathcal{T}_l^{C} \subseteq \mathcal{T}^{T}_l$, s.t. $|\mathcal{T}_l^T| = m_l$ where $m_l \sim \text{Uniform}(m_{min}, m_{max})$.
        }

        Concatenate to obtain the \hspace{30mm} \textbf{Context dataset}, $\mathbb{C} = \{\mathcal{T}_i^{C}\}_{i=1}^{M}$, \hspace{30mm} \textbf{Target dataset}, $\mathbb{T} = \{\mathcal{T}_i^{T}\}_{i=1}^{M}$.
    }

    \tcp{Model prediction and optimisation}
    \For{$j = 1$ \KwTo $N_{\text{sys}}$}{
        \For{$k = 1$ \KwTo $N_{x_0}$}{
            Sample ${task} \sim \text{Bernoulli}(\lambda)$.\;

            Obtain the new context,
            \eIf{${task} = 1$}{
                \tcp{Forecast}
                $\mathbb{C}_{\text{update}} = \{(t_{0,k}, \bm{x}_{0,k})\}$ where $(t_{0,k}, \bm{x}_{0,k}) \in \mathcal{T}_k$.\;
            }{
                \tcp{Interpolate}
                \If{$k \geq M$}{
                Randomly subsample $\mathcal{T}_k$ to get the context $\mathcal{T}_k^C$, and target $\mathcal{T}_k^T$ as in the training stage. \;
                }
                $\mathbb{C}_{\text{update}} = \{\mathcal{T}_k^C\}$
                
            }

            Obtain the new targets $\mathbb{T}_{\text{update}} = \{\mathcal{T}_k^T\}$
            
            Update context and target sets: \hspace{30mm}
            $\mathbb{C} \leftarrow \mathbb{C} \cup \mathbb{C} _{\text{update}}$, \hspace{30mm}
            $\mathbb{T} \leftarrow \mathbb{T} \cup \mathbb{T} _{\text{update}}$.\;

            \textbf{Compute} the variational posteriors. \;
            \textbf{Sample} the latent space using the variational posteriors. \;
            \textbf{Solve} the latent ODEs as given in \autoref{eq:ODEevolve}. \;
            \textbf{Decode} from latent space to obtain model predictions on the target space. \;
            Compute the \textbf{trajectory loss} $\mathcal{L}_{\text{ELBO},k}$ given in \autoref{eq:Liklihood_SANODEP}.\;
        }
    }

    Compute \textbf{mean loss}: $\mathcal{L}_{\text{ELBO}} = \frac{1}{N_{x_0}} \sum_k \mathcal{L}_{\text{ELBO},k}$.\;
    Update \textbf{model parameters}: $\theta \leftarrow \theta - \eta \nabla_\theta \mathcal{L}_{\text{ELBO}}$.\;
}
\end{algorithm*}

\section{Case Study: Penicillin Production} 
\label{sec:case_study}

Fed-batch systems combine an initial concentration of unreacted substrates and a feed of reactant into a reactor to dynamically generate a desired product at the end of the batch. Fed-batch reactor optimisation is a longstanding subject in the process systems engineering literature, with problems including feed profile \citep{cuthrell1989simultaneous}, initial condition \citep{zhou2009optimizing}, and stopping time \citep{patron2024economic,patron2024economically} optimisation. As a prototypical example, the production of penicillin in a fed-batch reactor has been used widely to demonstrate  control \citep{cuthrell1989simultaneous,lim1986computational} and monitoring principles \citep{shokry2018data} under varying initial conditions. In this case study, we consider the aforementioned optimisation of feed, initial condition, and stopping time for the penicillin reactor in a few-shot, i.e., few experiments, setting enabled by SANODEP. Crucially, SANODEP also enables the learning of parameter distributions to hedge against batch-to-batch fluctuations in exogenously varying stochastic reaction parameters. The mechanistic penicillin reactor model is described next.

\subsection{ODE System} \label{sec:case_study_ode_system}

A substrate feed with concentration $S_{F}(g/L)$ and flowrate $F(g/hr)$ is input into a reactor to initiate biomass growth resulting in penicillin production. The dynamical state of the system is modelled through the concentrations of biomass $B(g/L)$, substrate $S(g/L)$, and penicillin $P(g/L)$, as well as the reactor liquid volume $V(L)$. The reactor is seeded with an initial mixture defined by the initial conditions: $B_{0}$, $P_{0}$, $S_{0}$, and $V_{0}$. The process is then governed by the mass conservation laws~\citep{bajpai1980mechanistic} below: 

\begin{equation}\label{eq:1}
    \frac{dB}{dt} = \mu \left(B,S\right) B - \left(\frac{B}{S_{F} V}\right) F ;  \quad  B(t=t_0) = B_{0}
\end{equation}
\begin{equation}
\frac{dP}{dt} = \rho \left(S\right) B - K_{deg} P - \left(\frac{P}{S_{F} V}\right) F ;  \quad   P(t=t_0) = P_{0}
\end{equation}
\begin{equation}
\frac{dS}{dt} = - \mu \left(B,S\right) \left(\frac{B}{Y_{B/S}}\right) - \rho \left(S\right) \left(\frac{B}{Y_{P/S}}\right) - \gamma \left(S\right) B + \left(1 - \frac{S}{S_{F}}\right) \frac{F}{V} ; 
 S(t=t_0) = S_{0} 
 \end{equation}
\begin{equation}
\frac{dV}{dt} = \frac{F}{S_{F}} ;  \quad  V(t=t_0) = V_{0},
\end{equation}
where $\mu \left(B,S\right)(1/hr)$ is the biomass growth rate as described by a Contois model,  $\rho \left(S\right)(gP/gB/hr)$ is the penicillin production rate as described by a Monod model, and $\gamma \left(S\right)(gS/gB/hr)$ is also a Monod model used to describe the substrate requirement by the biomass. Specifically, these biochemical kinetics are given as:
\begin{equation}
\mu \left(B,S\right) = \mu_\mathrm{max} \left(\frac{S}{K_{B} B + S}\right),
\end{equation}
\begin{equation}
\rho \left(S\right) = \rho_\mathrm{max} \left(\frac{S}{K_{P} + S(1 + S/K_{in})}\right),
\end{equation}
\begin{equation}\label{eq:7}
\gamma \left(S\right) = m_{S} \left(\frac{S}{K_{m} + S}\right).
\end{equation}

\begin{figure}[h]
    \centering
    \includegraphics[width=.4\linewidth]{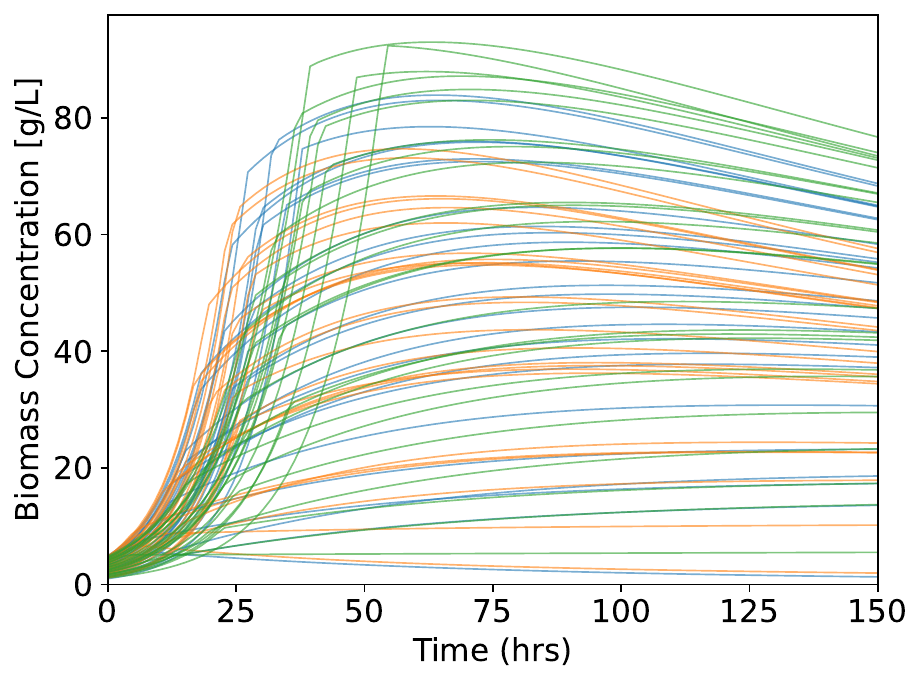} \includegraphics[width=.4\linewidth]{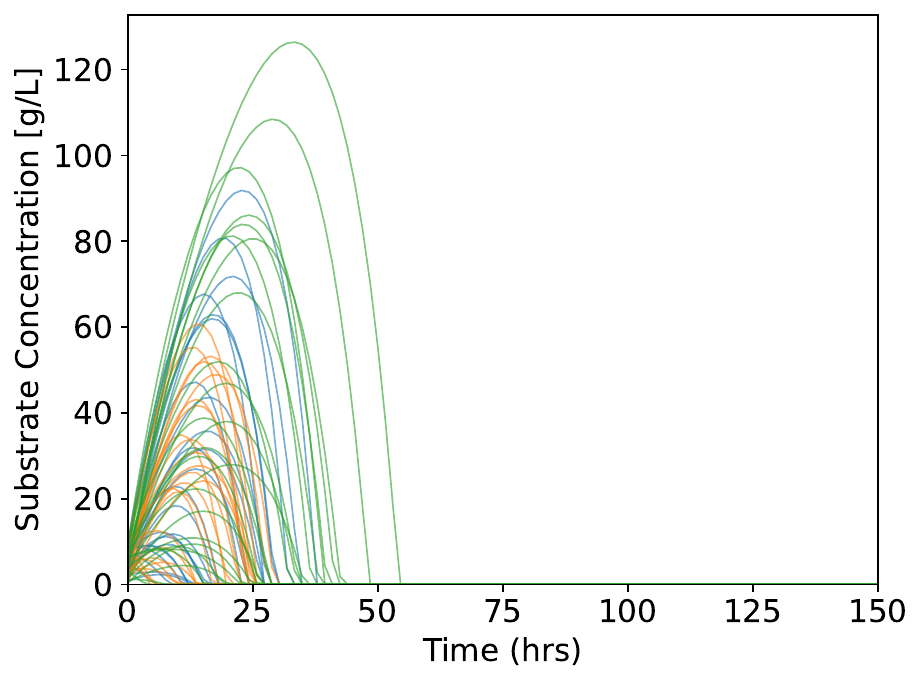}
    \\
    \includegraphics[width=.4\linewidth]{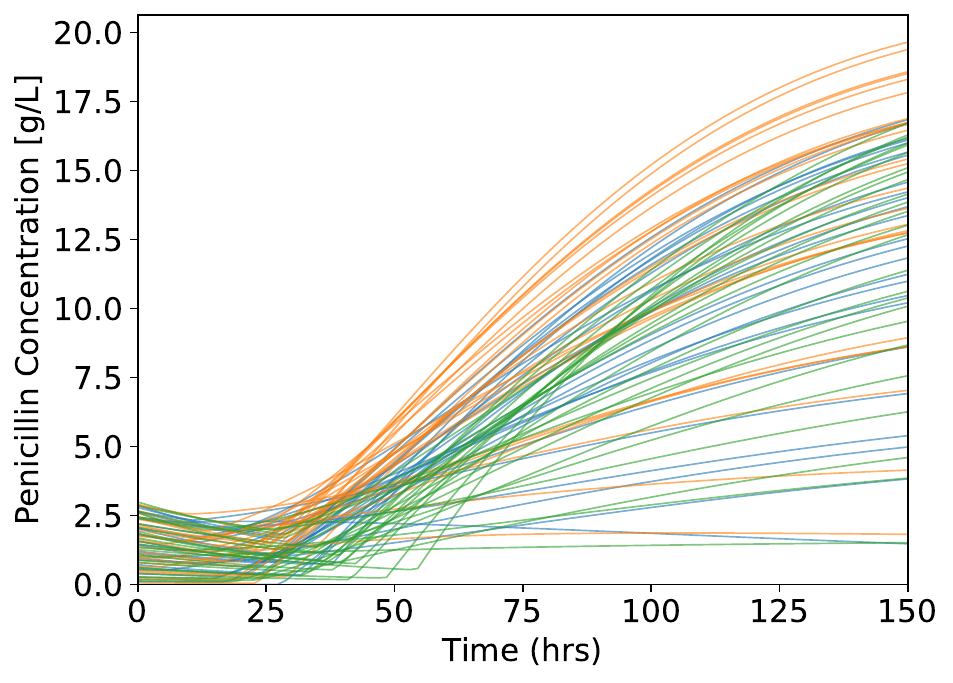} \includegraphics[width=.4\linewidth]{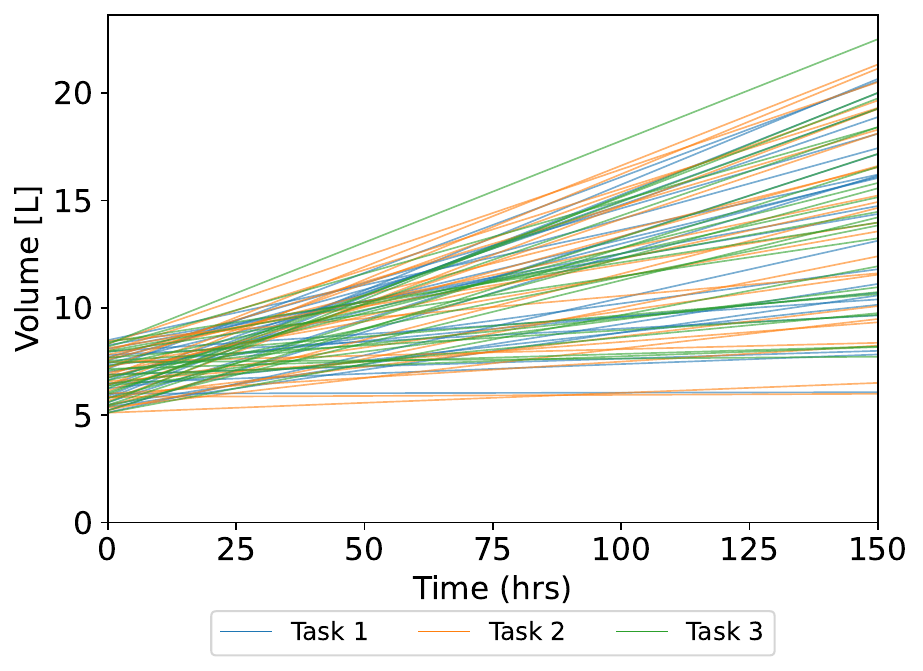}

    \caption{Example trajectories simulated using the ODE solver for randomly sampled initial conditions for each of the state variables against experiment time $[hrs]$;  Biomass Concentration $[g/L]$, Substrate Concentration $[g/L]$, Penicillin Concentration $[g/L]$,  and Volume $[V]$. Each of the three tasks presented is randomly sampled from the ranges derived from literature as in the training distribution given in \autoref{tab:stochastics_params}.}
    \label{fig:samplebatches}
\end{figure}

Definitions and values of fixed parameters are based on \citet{cuthrell1989simultaneous} and provided in Table \ref{tab:fixed_values}. The model parameters are assumed to vary from batch to batch and are modelled as uniform distributions of width $\sigma_{prior}$ around their nominal values (Table~\ref{tab:fixed_values}). Allowable ranges for the initial conditions are based on \citet{shokry2018data} and given in Table \ref{tab:init_cond_ranges}. ``Tasks'' $\bm{k}$, i.e., sets of values for the stochastic model parameters, are sampled across various initial conditions $\bm{x}_0$, which are also sampled from their respective  distributions. From these, a DAE solver is used to simulate the system \eqref{eq:1}--\eqref{eq:7}; a few example trajectories are shown in Figure~\ref{fig:samplebatches}. Each optimisation task seeks to optimise a selected objective given a set of the sampled stochastic parameters.

\subsection{Objective} \label{sec:case_study_objectives}

Productivity can be captured by a broad range of methods depending on the overall objective of the process. As a representative example, we take profit $g(\$)$ as the objective herein~\citep{profit_paper}:
\begin{equation}
    g = 2.5\times 10^{-2} \, P(t)\, V(t) - 168t - 8.5\times 10^{-4} \int_0^t F(z)\, dz.
    \label{eq:profit}
\end{equation}

\begin{table}[b]
  \caption{Fixed model parameters from \citep{cuthrell1989simultaneous}.}
  \label{params}
  \centering
  \begin{tabular}{lll}
    \toprule
    Symbol   & Description     & Value \\
    \midrule
    $S_{F}$ & Substrate concentration in feed  & $500 \; (g/L)$ \\
    $K_{deg}$  & Hydrolysis (degradation) constant & $0.01 \; (1/hr)$ \\
    $K_{in}$  & Substrate inhibition constant & $0.1 \; (gS/L)$ \\
    $Y_{B/S}$  & Biomass to substrate yield & $0.47 \; (gB/gS)$ \\
    $Y_{P/S}$  & Penicillin to substrate yield & $1.2 \; (gP/gS)$ \\
    \bottomrule
  \end{tabular}
  \label{tab:fixed_values}
\end{table}

\begin{table}[hb]
  \caption{Initial conditions and sample distributions from \citep{shokry2018data}.}
  \label{params}
  \centering
  \begin{tabular}{llll}
    \toprule
    Symbol     & Description     & Nominal value & Distribution\\
    & & (units) & \\
    \midrule
    $F$ & Substrate feed flowrate  & $25 \; (g/hr)$ & $\mathcal{U}_{[0,50]}$\\ 
    $B_0$ & Initial biomass concentration  & $1.5 \; (g/L)$ & $\mathcal{U}_{[1,5]}$\\ 
    $P_0$ & Initial penicillin concentration  & $0 \; (g/L)$ & $\mathcal{U}_{[0,3]}$ \\
    $S_0$ & Initial substrate concentration & $0 \; (g/L)$ & $\mathcal{U}_{[0,10]}$  \\
    $V_0$ & Initial liquid volume & $7 \; (L)$ & $\mathcal{U}_{[5,8.5]}$  \\
    \bottomrule
  \end{tabular}
  \label{tab:init_cond_ranges}
\end{table}

\begin{table}[htb]
\caption{Stochastic model parameters from  \citet{cuthrell1989simultaneous} and their nominal values.}
  \label{params}
  \centering
  \begin{tabular}{lll}
    \toprule
    Symbol     & Description     & Nominal value (Units) \\
    \midrule
    $K_{B}$ & Contois saturation constant for biomass production & $0.006 \; (gS/gB)$ \\
    $K_{P}$  & Saturation constant for substrate consumption & $0.0001 \; (gS/L)$ \\
    $K_{m}$  & Monod saturation constant for substrate maintenance & $0.0001 \; (gS/L)$ \\
    $\mu_{max}$  & Maximum biomass growth rate & $0.11 \; (1/hr)$ \\
    $\rho_{max}$  & Maximum penicillin production rate & $0.0055 \;p (gP/gB/hr)$ \\
    $m_{S}$ & Maintenance requirement of substrate by biomass  & $0.029 \; (gS/gB/hr)$ \\
    \bottomrule
  \end{tabular}
  \label{tab:stochastics_params}
\end{table}

\section{Results} \label{sec:results} 

This section presents comprehensive simulated results demonstrating the impact of meta-learned surrogate models on optimisation performance within the outlined biochemical domain. In our analysis, we compare across the below BayesOpt strategies: 

\begin{itemize}
    \item[(i)] \textbf{SANODEP}, which uses SANODEP as the surrogate model, outlined in \autoref{sec:model_model_structure}, pre-trained using the bespoke training procedure, and using the adaptable acquisition function, both outlined in \autoref{sec:model_bayes_op}. The standard baseline prior, $\bm{k} \sim \pi_{k} ; \; \mathcal{U}\left[(1-\Delta_k)\bm{k}_\mathrm{nom}, (1-\Delta_k) \bm{k}_\mathrm{nom}\right] $, is used to generate the \textbf{meta-training data} where the industry values, $\textbf{k}_\mathrm{nom}$, are outlined in \autoref{sec:case_study}. The window is selected to be $\Delta_k = 0.05$.
    \item[(ii)] \textbf{GP-Standard} uses a standard zero-mean GP with an RBF kernel for the surrogate model alongside the Expected Improvement (EI) acquisition function introduced in \autoref{sec:background}. This provides a problem-agnostic baseline standard to compare against. Note that GP-BO takes only a single sample observation at the time $t^*$ which optimises the acquisition function, unlike SANODEP which has an adaptable acquisition function allowing for multiple samples from the same trajectory. 
    \item[(iii)] \textbf{GP–Exp}, which augments GP-BO observations with additional intermediate trajectory samples uniformly sampled up to the selected stopping point. In the context of this case study, intermediate samples are inexpensive in comparison to the cost of sampling a new, unseen trajectory and so are provided to the model at no added cost. 
\end{itemize}

A key factor determining the efficacy of meta-learning approaches is the relevance of the dataset used to train the base model. 
Here, we train the base model on samples from the true ODE system considering parametric variations: $\bm{k} \sim \pi_{k} ; \; \mathcal{U}\left[(1-\Delta_k)\bm{k}_\mathrm{nom}, (1-\Delta_k) \bm{k}_\mathrm{nom}\right] $. 
While this is a strong assumption, we study the effect of the meta-learning dataset quality by testing on tasks outside of $\bm{k}$, i.e., tasks the base model has never seen. 
Specifically, we assess the performance of SANODEP using on- and off-task testing data generated by introducing an offset to the mean of the prior distribution, $\delta$, over the nominal task values and narrowing the window, $\Delta_{\delta}$. 
The resulting testing distribution is given,
\begin{equation}
    \bm{k}_\mathrm{test} \sim \pi_{\delta} ; \; \delta\bm{k}_\mathrm{nom} + \mathcal{U}\left[(1-\Delta_{\delta})\bm{k}_\mathrm{nom}, (1 + \Delta_{\delta}) \bm{k}_\mathrm{nom}\right],
\end{equation} 
where the similarity between $\bm{k}$ and $\bm{k}_{test}$ effectively determine the relevance of the meta-learning dataset, including when $\bm{k}$ and $\bm{k}_{test}$ have zero overlap. Future work may study the effect of structural mismatches in the model. 
The regions considered in our analysis are presented in table \autoref{tab:testing_prior}.

\begin{table*}[h!] 
  \caption{Parameters of the training, $\pi_k$, and testing, $\pi_{\delta}$ distributions presented in this section. Positive offsets locate the distribution at, or off, the upper edge. Negative offsets locate the distribution at, or off, the lower edge.}
  \centering
  \begin{tabular}{lllll}
    \toprule
    Distribution     & Type   & Description  & Offset $\delta$ & Window $\Delta$ \\
    \midrule
    On-Task; $\pi_k$ & Training & & $0$ & $0.05$ \\
    Very Off Task & Testing & & $\pm 0.5$ & $0.01$ \\
    Slightly Off Task & Testing & & $\pm 0.06$ & $0.01$ \\
    Almost Off Task & Testing & & $\pm 0.04$ & $0.01$ \\
    On Task & Testing & & $0.00$ & $0.01$ \\
    \bottomrule
    \end{tabular}
\label{tab:testing_prior}
\end{table*}

\subsection{Modelling Results}

\begin{figure}[]
    \centering
    
    \raisebox{-0.5\height}{\includegraphics[width=0.35\linewidth]{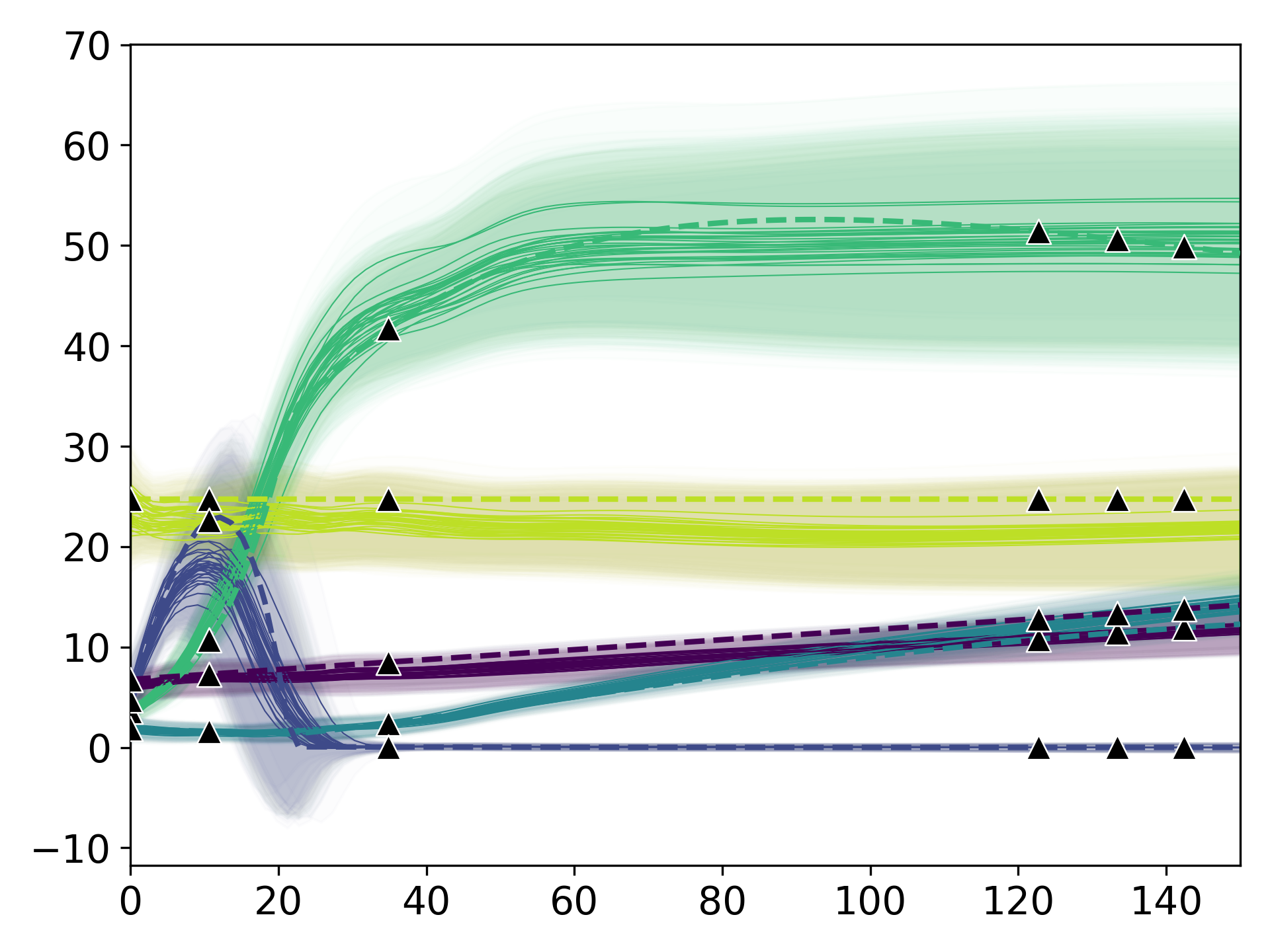} }
    \raisebox{-0.5\height}{\includegraphics[width=0.35\linewidth]{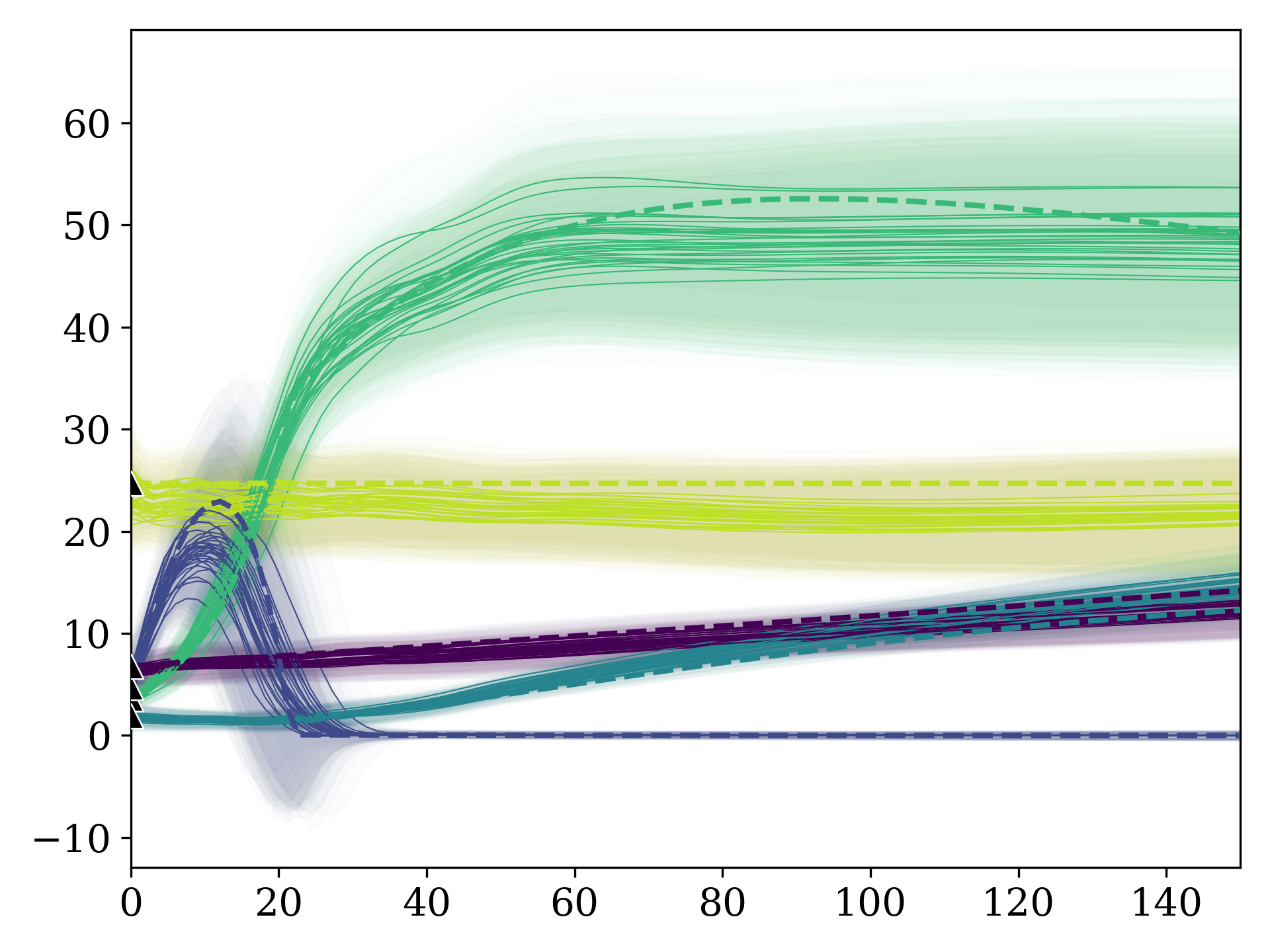} }
    \raisebox{-0.5\height}{\includegraphics[width=0.2\linewidth]{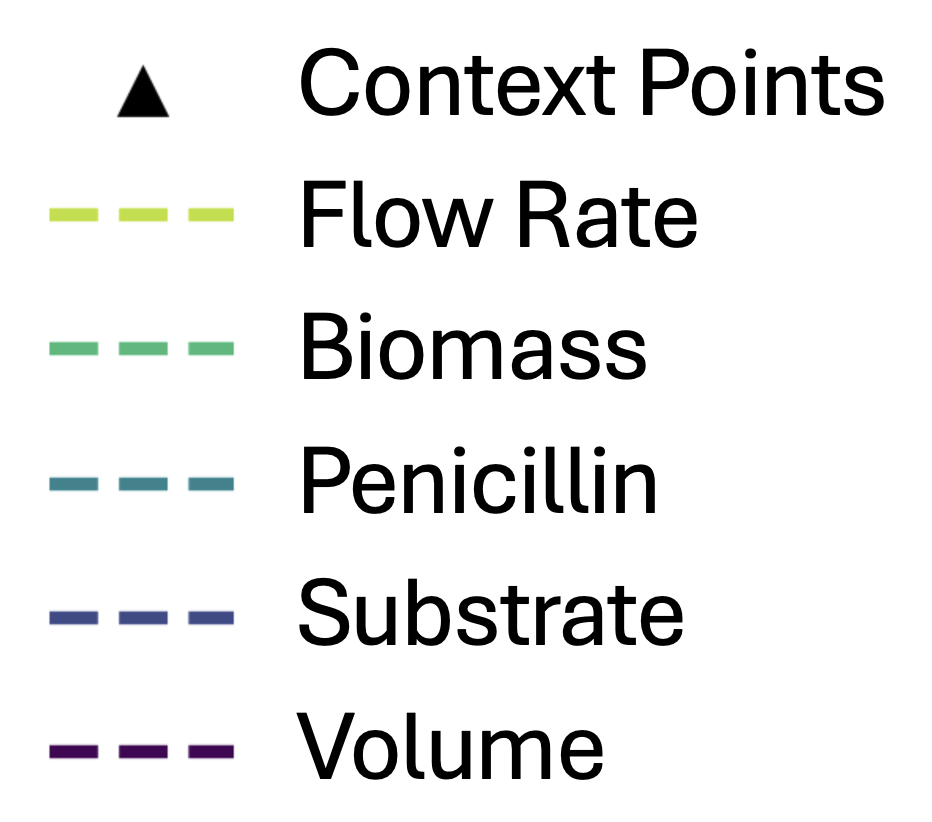}}
    \caption{Example trajectories drawn from a single task from the On-Task distribution using a pretrained SANODEP model. On the right, the model is in forecast mode, and on the left, the model is in interpolation mode. Observations are indicated with a triangle.}
    \label{fig:interpolation_forecast_sanodep}
\end{figure}

SANODEP supports two prediction modes: forecasting and interpolation. Forecasting corresponds to predicting an entire state trajectory given only the initial conditions, whilst interpolation incorporates some intermediate observation(s) taken during the trajectory. In both modes, predictions may be made for a task for which either none or some prior trajectories have been observed. BayesOpt harnesses both of these modes and makes predictions on trajectories within one task; i.e. using the same stochastic model parameters sampled at the start of the BayesOpt run. Examples of these trajectories for both modes are given in \autoref{fig:interpolation_forecast_sanodep}.

\begin{figure}[h!]
    \centering
    \includegraphics[width=0.5\linewidth]{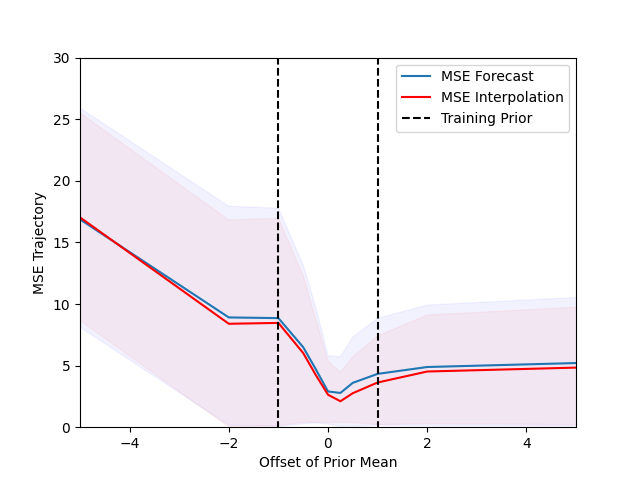}
    \caption{Trajectory-wise MSE of samples drawn from SANODEP across tasks with varying offsets from the nominal prior in the penicillin batch-reactor setting. MSE calculation is done on twenty trajectories drawing from the MSE testing prior. Descriptions of the regions plotted vertically are given in \autoref{tab:testing_prior}. MSE is highest for negative offsets, decreases toward the nominal task, and is lowest at the centre of the prior. }
    \label{fig:mse_prior_offset}
\end{figure}

On-task and off-task performance are both measured using the Mean Squared Error (MSE) between predictions and the true trajectory for tasks sampled using the distributions outlined. Tasks are sampled using the standard testing window $\Delta_{\delta} = 0.01$ varying the offset $\delta$. In \autoref{fig:mse_prior_offset} we present the results of this analysis.

These results indicate that SANODEP has been successfully trained on the meta-training dataset generated from the specified distributions. On-task performance, measured via MSE for samples generated within the nominal prior, is strongest near the centre of the task. This is expected, as central samples have more neighbouring points within the meta-training dataset that are similar, whereas samples near the edges of the distribution are sparser, making central points effectively overrepresented during training. The problem setting additionally exhibits asymmetry, with negative offsets producing non-physical and unrealistic behaviour which can diverge. These findings demonstrate that a narrow training prior can capture system behaviour across a broad range of tasks, providing robust predictions for use in the optimisation of a range of tasks both on- and off-distribution. 

\subsection{Optimisation Results -- Base Case}

In this section we present the results of Bayesian Optimisation using the models outlined in the section above. Each iteration of the BayesOpt loop involves the updating of the surrogate, the evaluation of the acquisition, and the return of a sample from the black-box function, i.e., an experiment. For each full run of BayesOpt, a single task is selected from the indicated testing distribution, $\bm{k}_\mathrm{test}$ and is used as the ground-truth, black-box system. 
The number of experiments is dictated by the \textbf{sampling budget} which determines how many samples of the black-box function may be drawn for the duration of the run. In this work, we take an experiment, or black-box sample, to be a single trajectory, such that the sampling budget is the number of initial conditions $\bm{x}_0$ and stopping point $t_\mathrm{stop}$ pairs to be queried. Within a trajectory, model-dependent budgets are set for the number of additional intermediate samples $\bm{x}_t$, which may be taken at times $t_i$. 

For SANODEP, the meta-learning procedure is pre-computed offline using the meta-training prior before BayesOpt begins. The model is subsequently updated online using task-specific, in-context black-box observations collected during the BayesOpt loop. In contrast, GP-based surrogate models do not rely on offline meta-training and instead require an initial dataset drawn directly from the black-box function to initialise the model, i.e., fit the GP hyperparameters. In the analysis presented here, this initial dataset is generated using Latin Hypercube Sampling (LHS) over the input space defined by the initial conditions and time $(\bm{x}_0, t)$. It is important to note that, SANODEP provides trajectory predictions $\bm{x}_t$ which are used to calculate the objective $\textbf{g}(\bm{x}_t)$ where as GP-based methods model $\textbf{g}(\bm{x}_t)$ directly.

Direct cross-task comparisons are made difficult by the differing maxima that may be achieved across different tasks. In this analysis, we normalise the performance on each task using the task-specific maximum, found by taking the maximum observed objective value from a single run of GP-BO on the task with a large sample budget of 50. We then normalise the objective, defined in \autoref{sec:case_study_objectives}, by this task maximum to allow for cross-task comparison in the analysis presented in this section. 

The pre-training phase of SANODEP allows it to capture shared structure across tasks and operate effectively in low-data regimes. However, we note the expressive latent ODE structure and amortised inference required by SANODEP result in a substantially higher computational cost per optimisation step compared to simpler surrogate models. Once initialised, GP inference and acquisition function optimisation are comparatively inexpensive. Similarly, random search incurs negligible computational overhead. 
Nevertheless, given its sample efficiency, the sample budget is set to $10$ black-box evaluations for SANODEP and doubled to $20$ for the GP-based methods and the random search. Practically, we found this setting allows for all methods to reach close-to-optimal performance, roughly $\textbf{g}_\mathrm{norm} > 0.85$, where $\textbf{g}_\mathrm{norm}$ is the mean of the normalised objective maxima across five random repetitions.  

\begin{figure}[h!] 
    \centering
    \includegraphics[width=0.5\linewidth]{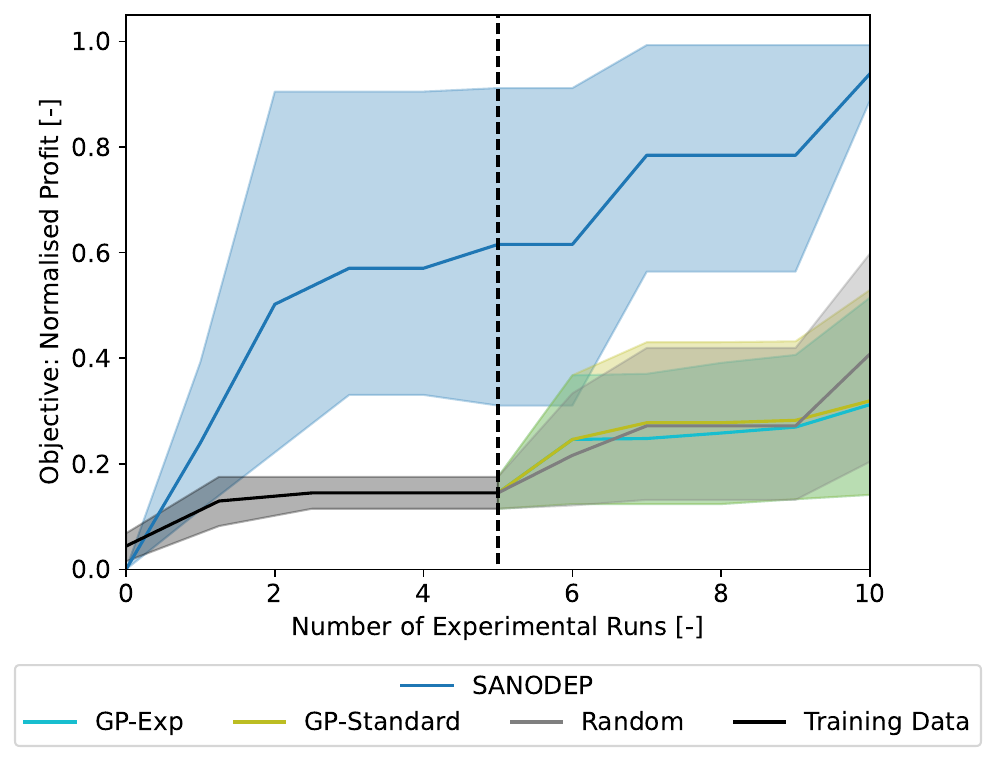}
    \caption{On-task optimisation performance for a single task sampled from the centre of the prior, measured against the number of experimental runs. SANODEP performs strongly in the low-data regime when compared to the traditional methods.}
    \label{fig:centre_task}
\end{figure}

\begin{figure}[h!]
    \centering
    \includegraphics[width=0.5\linewidth]{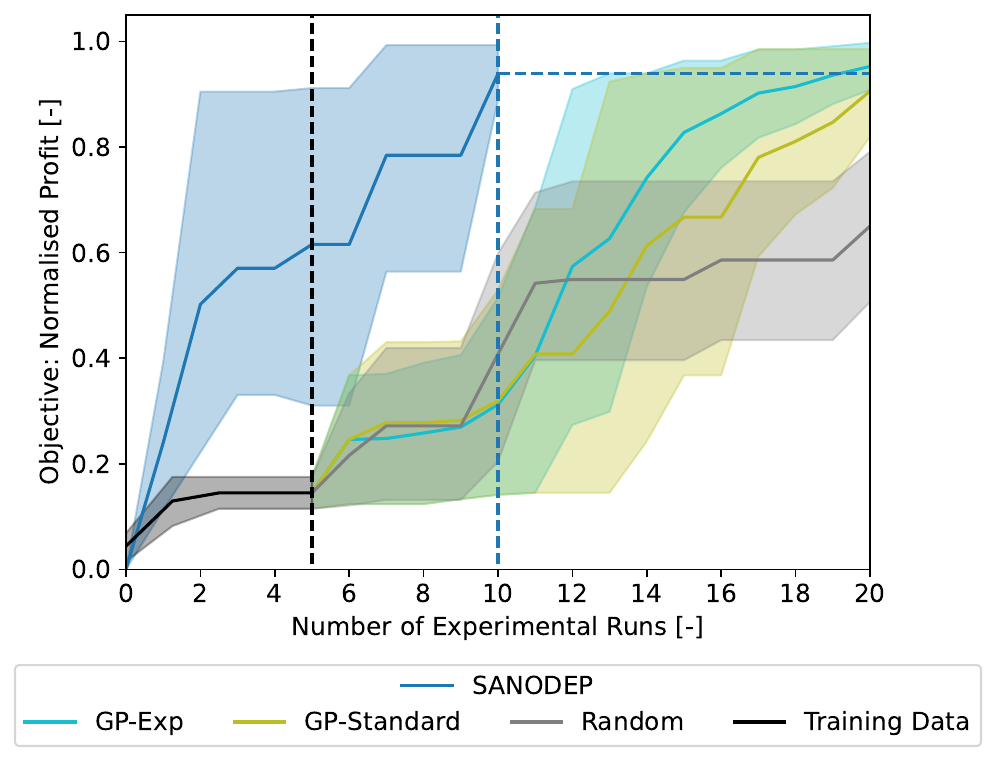}
    \caption{On-task optimisation performance in the ``infinitum'' presented in \autoref{fig:centre_task}. Performance at this boundary is close to optimum. It was experimentally observed that extension beyond these horizons is computationally expensive and results in marginal gain in the normalised regime. In the 'infinitum', GP perform on par with SANODEP.}
    \label{fig:centre_task_infinitum}
\end{figure}

\begin{figure}
    \centering
    \includegraphics[width=0.5\linewidth,trim=0 0 0 0.65cm, clip]{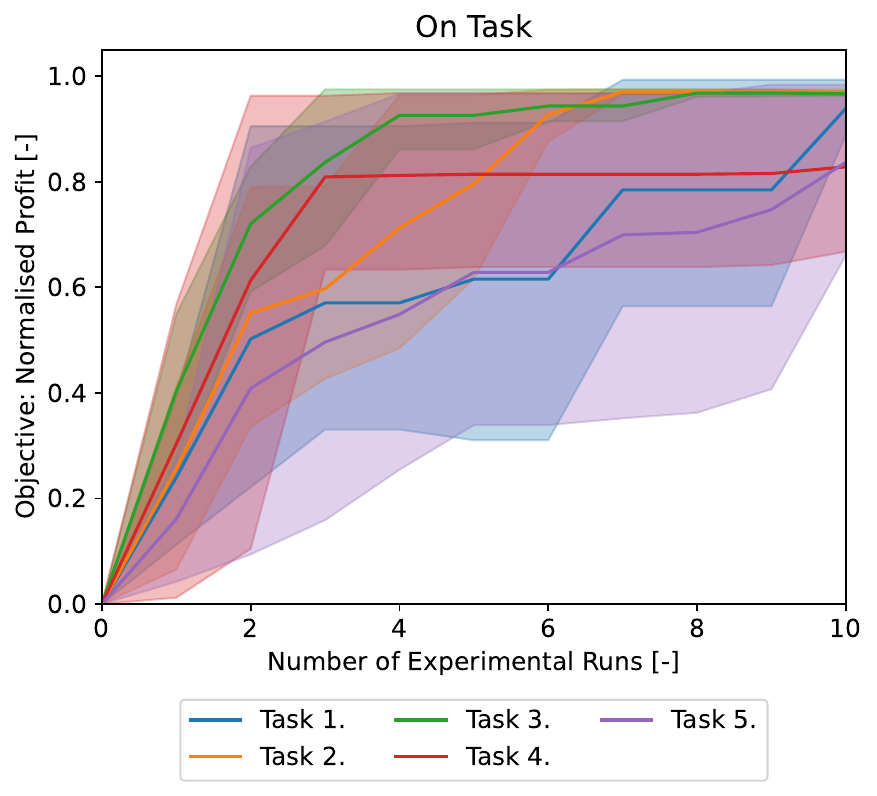}
    \caption{On-task SANODEP optimisation performance for $5$ tasks sampled from the On-Task distribution. The figure presents normalised profit against number of experimental runs. Performance across tasks varies significantly in the deviation within task but converges as the number of observed trajectories increases.}
    \label{fig:all_centre_tasks}
\end{figure}

We define the low-data regime to be $10$ trajectories, and compare performance of the methods outlined using a set of $5$ tasks drawn from the On-Task distribution defined in \autoref{tab:testing_prior}. Analysis for one of these tasks is presented in \autoref{fig:centre_task} for the low data regime alongside \autoref{fig:centre_task_infinitum} which presents further analysis on the same task for the complete sampling budget allocated to each method. For the full cohort of given tasks, SANODEP-specific analysis is presented in \autoref{fig:all_centre_tasks}. 
For this task, which lies squarely in the center of the parameter range used to pre-train the SANODEP model, SANODEP exhibits promising performance as a few-shot learner. 
In particular, BayesOpt with SANODEP appears to reliably optimise the system after only 10 experimental runs, converging to a similar point as the GP-based methods after the latter have queried 20 experiments. 

\begin{figure}[h!]
    \centering
    \includegraphics[width=0.49\linewidth]{Figures/_centre.pdf}
    \includegraphics[width=0.49\linewidth,trim=0 0 0 0.65cm,clip]{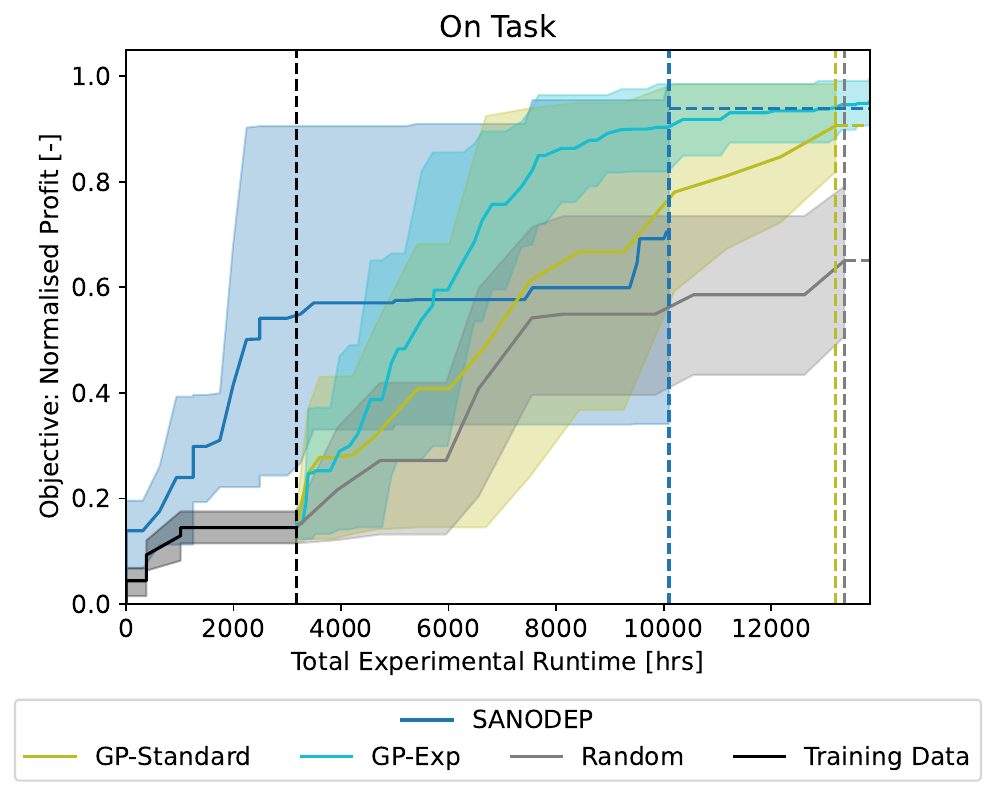}
    
    \caption{Plots show two common cost metrics in Bayesian Optimisation: (left) the number of state observations and (right) the total experimental runtime, both plotted against the optimisation objective, profit. These results include the offline, initial cost of collecting training data required for GP-based methods. SANODEP underperforms in this setting, as it is evaluated initially without access to any task-specific data. Depending on the application, either start-up reactor costs or experimental runtime may dominate the overall optimisation cost for which the optimisation strategy could be adapted to.}
    \label{fig:centre_time}
\end{figure}

Here, we note that the cost metric used to evaluate  the sampling budget is key to selecting a method. For the presented case study (and many batch reaction settings in general), the prohibitive cost is the number of batches run, i.e., separate initialisations of each trajectory. The length of time that each trajectory runs for is not as significant, and, subsequently, methods are not rewarded for choosing to end a trajectory early. Likewise, we assume the cost of taking measurement during a batch process to be relatively lower than starting a new batch. Analyses across two different cost metrics, the number of trajectories and the total experimental runtime, are presented in \autoref{fig:centre_time}.

\subsection{Optimisation Results -- Off-Task Setting}

In practical applications, while the historical data or prior modelling knowledge may determine a good meta-learning dataset, they may not always encompass the exact dynamics of the observed system. 
In other words, we may envision that the true values of the uncertain parameters (\autoref{tab:stochastics_params}) may lie outside the values used in the meta-learning prior. We note that structural mismatches are not considered here, but may be effectively approximated by a large mismatch in parameter values. Moreover, the computational budget required to train SANODEP effectively grows with the complexity of the training distribution. Given this motivation, we now evaluate model performance on tasks that extend beyond the training prior, providing insight into robustness and generalisation under realistic experimental conditions.

We define the edges of the training distribution using the training distribution $\bm{k}$. The upper edge is defined as the maximum parameter-wise values found in $\pi_{\bm{k}}$, and the lower edge is defined as the minimum parameter-wise values found in $\pi_{\bm{k}}$. We define three regions about these edges; Almost Off Task, Slightly Off Task, and Very Off Task. The full definitions of these are given in \autoref{tab:testing_prior}. The Almost Off Task distribution is fully contained within the training prior but located at the very edge of the boundary. The Slightly Off Task distribution lies completely outside of the training prior but is closely located to the edge. Trajectory samples from Slightly Off Task may closely resemble those from Almost Off Task, but notably were not present in the training prior. The Very Off Task distribution is significantly further from the edge and comprises tasks that are very different from the meta-learning dataset. 

\begin{figure}[h!]
    \centering
    \includegraphics[width=0.48\linewidth]{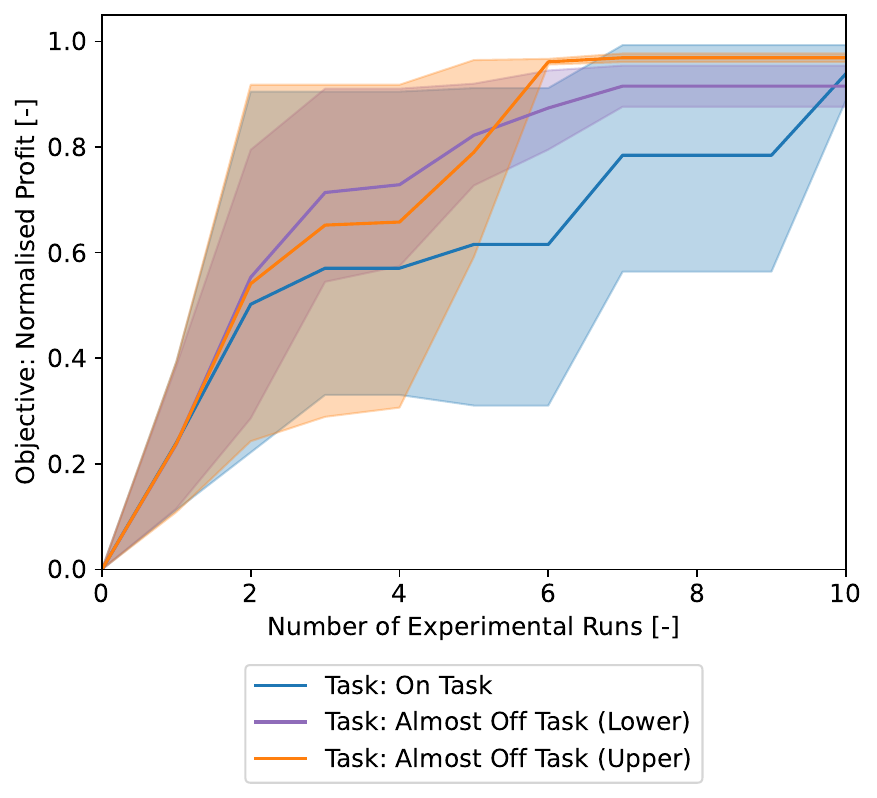}
    \includegraphics[width=0.48\linewidth]{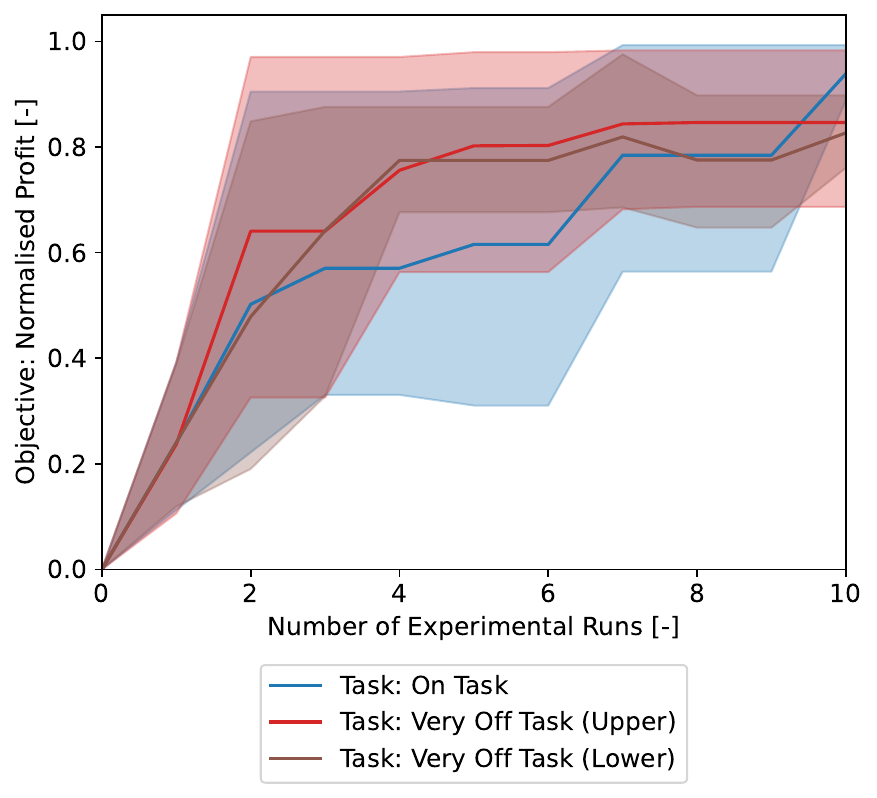} 
    \caption{Comparison of SANODEP performance on individual tasks sampled from the distributions distributions both on- and off-task, shown left and right respectively. }
    \label{fig:off_vs_on_task}
\end{figure}

\begin{figure*}[h!]
    \centering
    \includegraphics[width=0.4\linewidth]{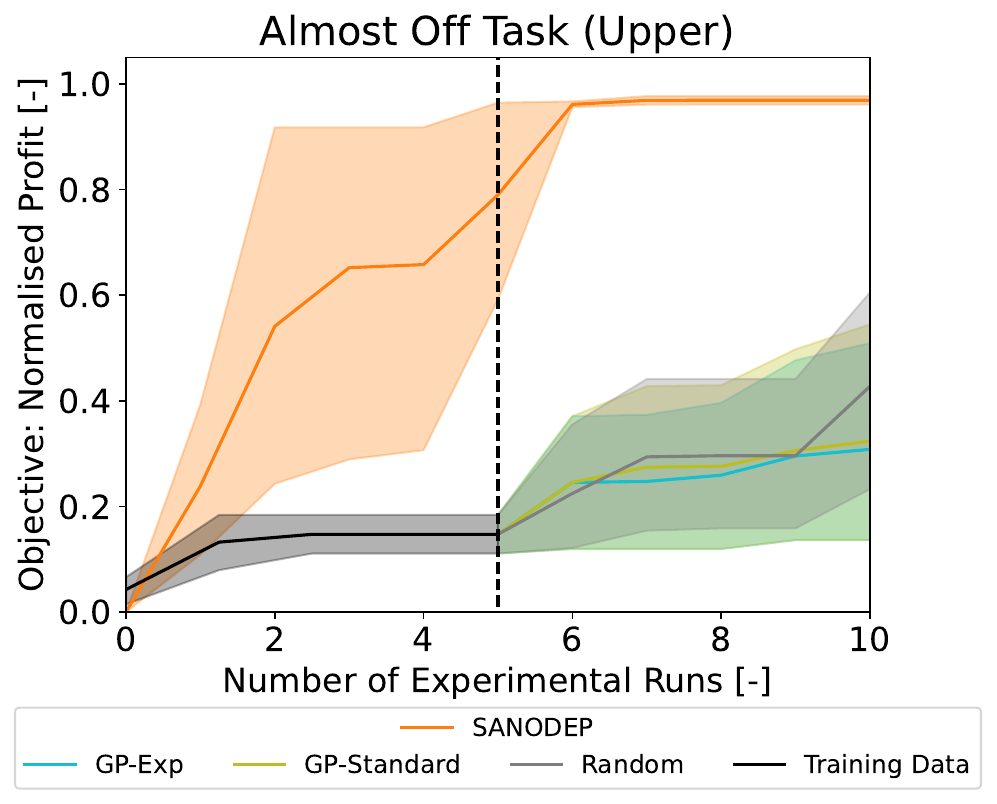} 
    \includegraphics[width=0.4\linewidth]{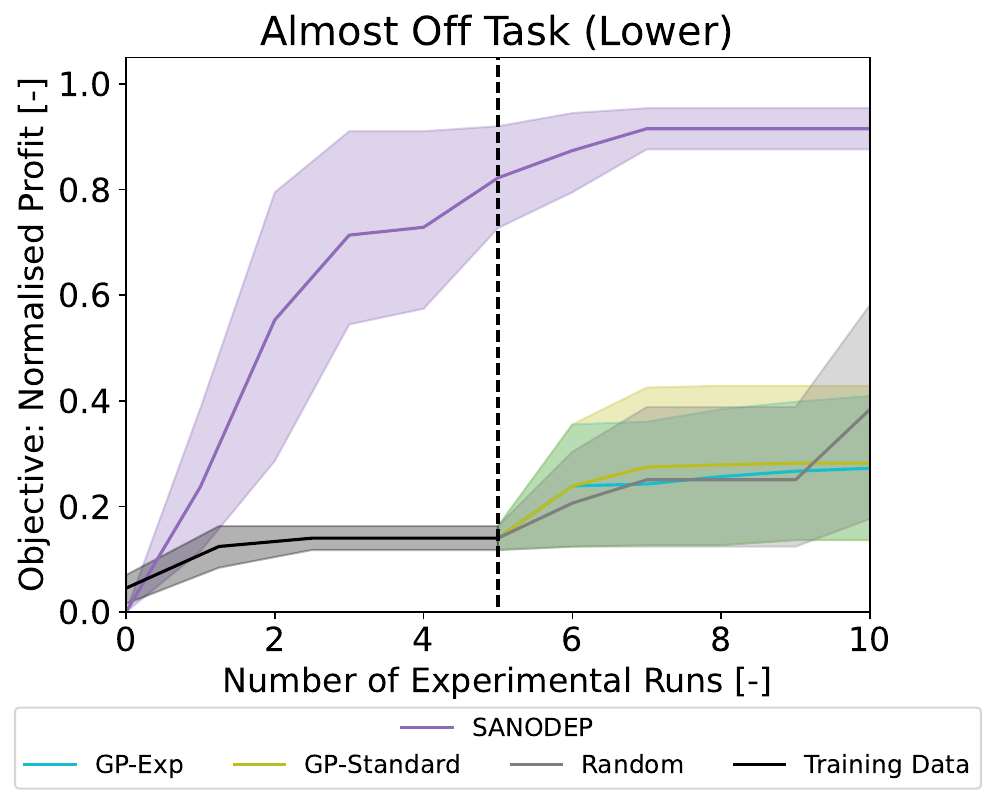} \\
    \includegraphics[width=0.4\linewidth]{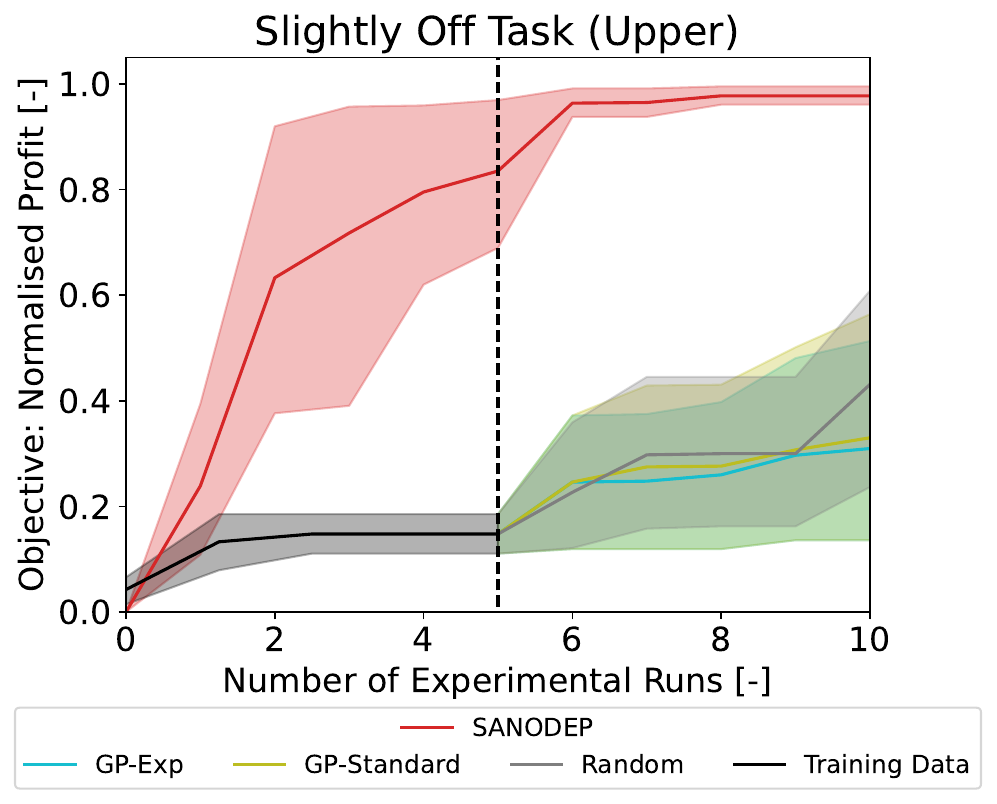} 
    \includegraphics[width=0.4\linewidth]{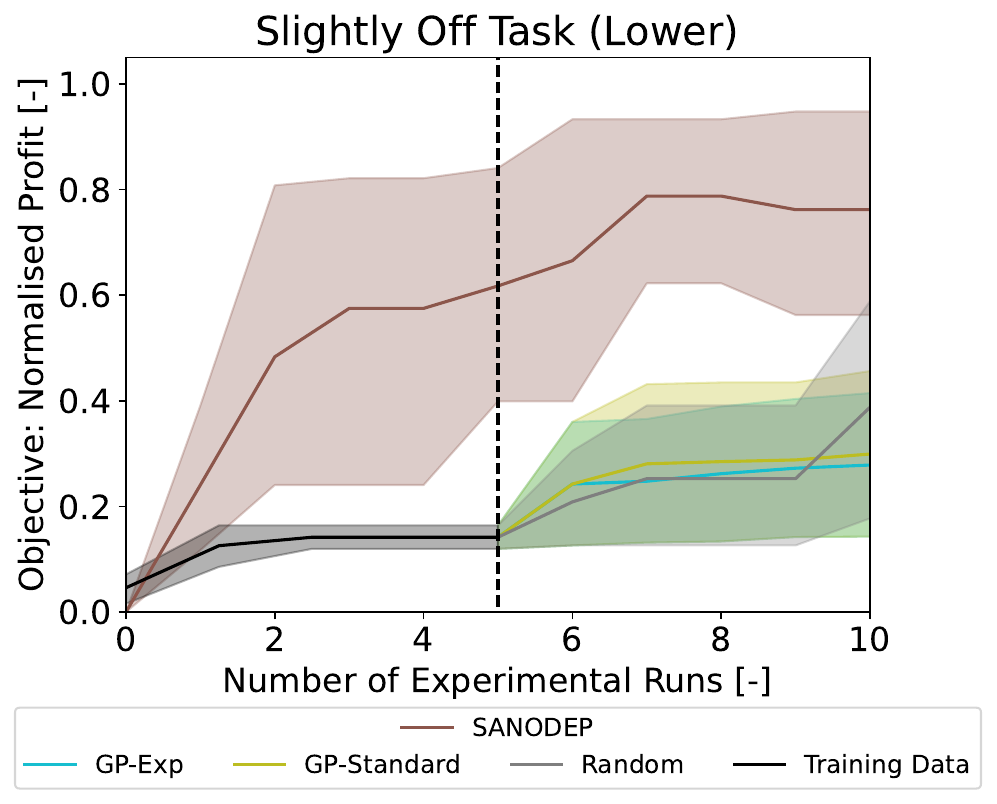} \\
    \includegraphics[width=0.4\linewidth]{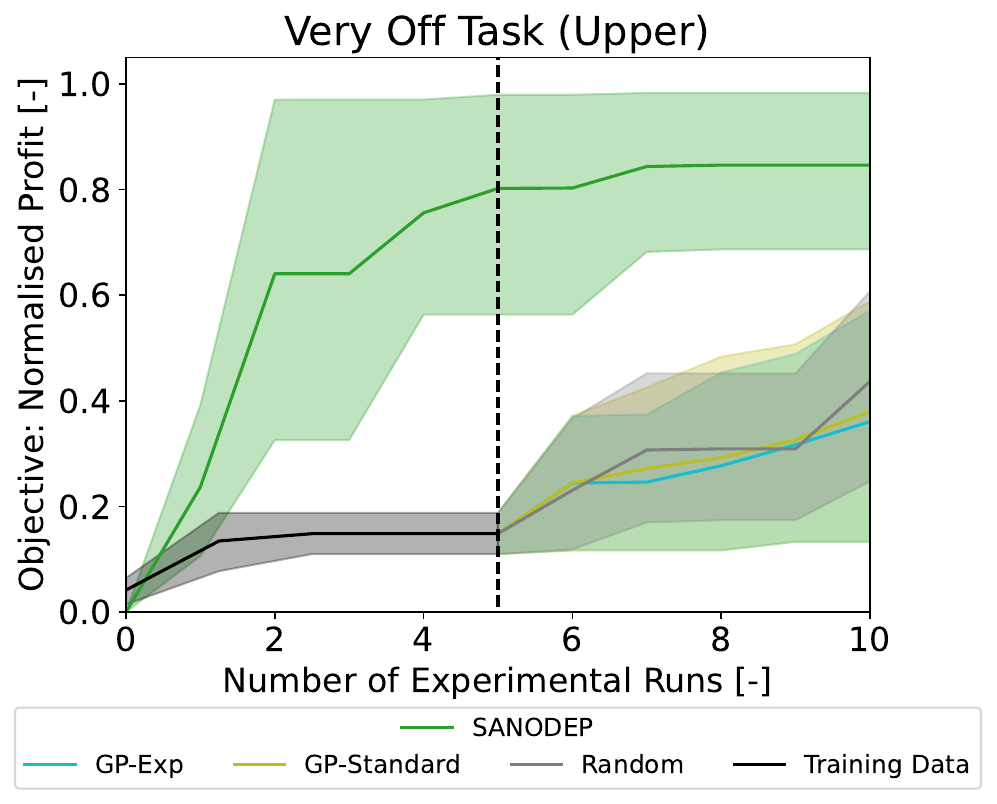} 
    \includegraphics[width=0.4\linewidth]{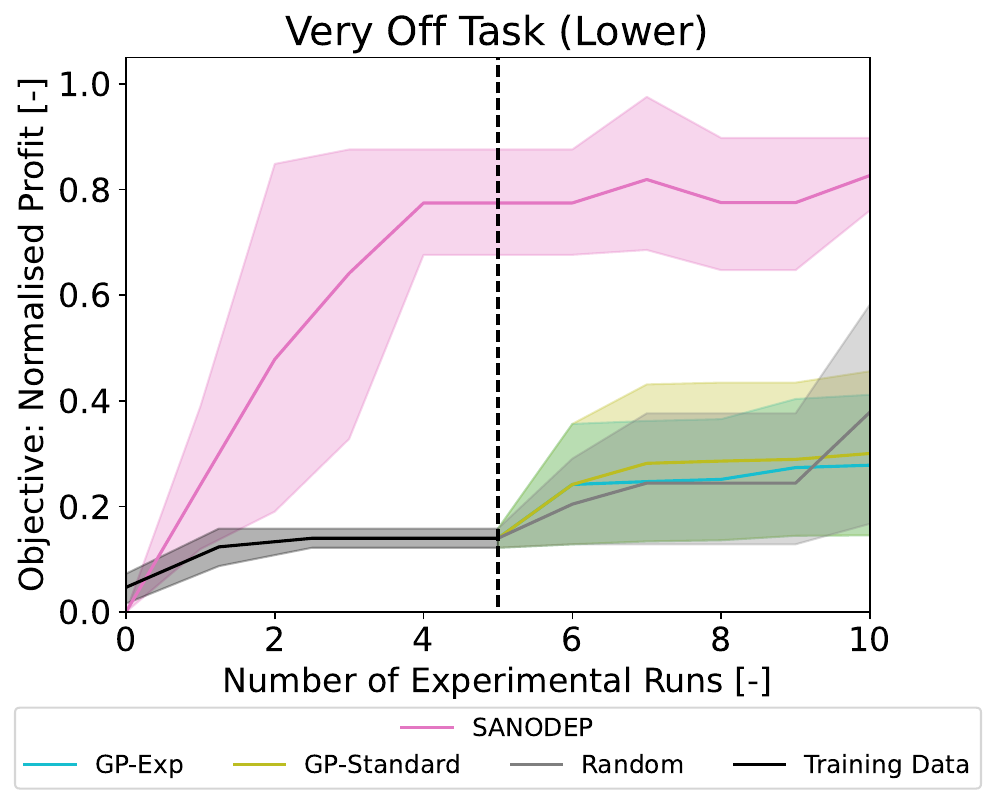} 
    \caption{Comparison of surrogate model performance across on- and off-task regimes. Each plot presents analysis for individual tasks which are sampled from distributions outlined in \autoref{tab:testing_prior}. Upper edge tasks are shown in left-hand column, and lower edge ones on right. }
    \label{fig:off_vs_on_task_gp}
\end{figure*}

Analysis for tasks sampled from these regimes is presented in \autoref{fig:off_vs_on_task} and compared against the GP-based methods and Random search in \autoref{fig:off_vs_on_task_gp}. 
\autoref{fig:off_vs_on_task} shows that SANODEP performance on tasks in the centre of the distribution is worse than on tasks located close to the edge. The model converges on the task-optima slower for the central task as opposed to those on the edge. This phenomenon may suggest that tasks closer to the upper/lower bound of parameter values are more `distinct.' In other words, the surrogate model can more quickly identify which trajectories in the base dataset are similar to the target system, using this knowledge to quickly optimise the given task. We note that performance near the lower edge of parameter values tends to be worse in comparison to tasks near the upper edge. This may occur because lower parameter values correspond to slower reaction rates. 

Finally,~\autoref{fig:off_vs_on_task_gp} shows that, while SANODEP tends to perform worse in very off-task settings, it still significantly outperforms the GP-Based methods in all regimes. Specifically, SANODEP-based BayesOpt shows decreased relative performance for both Very Off Tasks and the Lower Slightly Off Task, in terms of both the quality of the optima found, and the reliability of the convergence (higher variance in objective found after 10 experiments). 
We observe that the additional prior knowledge embedded in the meta-learning step appears structural: despite having not observed the task systems in the very off task setting, SANODEP can still quickly optimise the batch recipe in only a few experiments. 
This suggests that embedding the general mechanistic knowledge, e.g., which/how variables are related dynamically, provide useful information to the BayesOpt procedure, even when the actual data/models are inaccurate. 

\section{Conclusions} \label{sec:conclusion}

Fed-batch recipes are difficult to optimise experimentally when experiments are expensive and unmeasurable fluctuations in reaction conditions occur. Herein, we propose the use of SANODEP as an alternative to GPs for the few-shot Bayesian optimisation of these fed batch processes. SANODEP deploys meta-learning to improve generalisability to the unmeasurable fluctuations, as well as a dynamic structure to incorporate intermediate batch measurements into the model. We benchmark SANODEP against GPs in the Bayesian optimisation of a penicillin production process with stochastic reaction parameters. 

In the context of fed-batch chemical processes, SANODEP demonstrates strong performance in the low-data regime. It achieves this by adapting to newly observed batch trajectories through leveraging meta-learned priors derived from simulated training data. Analysis of mean squared error indicates that predictive accuracy of SANODEP as a surrogate model may decrease when task data lies outside of the training regime, which in turn can reduce BayesOpt performance, highlighting areas for further investigation in experimental design. However, when predictive accuracy is high which may occur off-task, SANODEP exhibits superior relative performance, and across all regimes, it outperforms GPs in the low-data scenario. GP-based Bayesian optimisation initially lags in few-shot settings but benefits from consistency guarantees, converging reliably to the process optimum as more batch observations are collected. This illustrates a trade-off between short-term adaptability and long-term convergence: SANODEP excels in scarce-data regimes, while GP-based methods catch up as data accumulate. Overall, these results emphasise the potential of meta-learning frameworks like SANODEP to accelerate efficient, reliable optimisation of batch chemical processes, enabling faster identification of optimal recipes and operating conditions in high-value manufacturing. Further, SANODEP could be used as an initial strategy to optimise conditions in the low-data regime before switching to GPs in the high-data regime.

\section*{Acknowledgments}
The authors gratefully acknowledge support from BASF SE, Engineering
and Physical Sciences Research Council [grants  EP/X025292/1, EP/Y028775/1, and EP/S023151/1] (RM, CT, JQ, BL), a BASF/RAEng Research Chair in Data-Driven Optimisation (RM), and a BASF/RAEng Senior Research Fellowship (CT). 
JQ acknowledges Research England’s Expanding Excellence in England (E3) fund awarded to MARS (Mathematics for AI in Real-world Systems) at Lancaster University. 
RM holds concurrent appointments as a Professor at Imperial College London and as an Amazon Scholar. This paper describes work performed at Imperial College London and is not associated with Amazon. 

\bibliography{rsc_full}

\begin{thebibliography}{60}
\providecommand{\natexlab}[1]{#1}
\providecommand{\url}[1]{\texttt{#1}}
\expandafter\ifx\csname urlstyle\endcsname\relax
  \providecommand{\doi}[1]{doi: #1}\else
  \providecommand{\doi}{doi: \begingroup \urlstyle{rm}\Url}\fi

\bibitem[Abdollahi and Dubljevic(2012)]{abdollahi2012lipid}
Javad Abdollahi and Stevan Dubljevic.
\newblock Lipid production optimization and optimal control of heterotrophic
  microalgae fed-batch bioreactor.
\newblock \emph{Chemical Engineering Science}, 84:\penalty0 619--627, 2012.

\bibitem[Ashraf and Rao(2022)]{ashraf2022multiobjective}
{Abdul Basith} Ashraf and {Chinta Sankar} Rao.
\newblock Multiobjective temperature trajectory optimization for unseeded batch
  cooling crystallization of aspirin.
\newblock \emph{Computers \& Chemical Engineering}, 160:\penalty0 107704, 2022.

\bibitem[Bajpai and Reuss(1980)]{bajpai1980mechanistic}
RK~Bajpai and M~Reuss.
\newblock A mechanistic model for penicillin production.
\newblock \emph{Journal of Chemical Technology \& Biotechnology}, 30\penalty0
  (1):\penalty0 332--344, 1980.

\bibitem[Barton et~al.(2021)Barton, Duran-Villalobos, and
  Lennox]{barton2021multivariate}
Maxwell Barton, Carlos~A Duran-Villalobos, and Barry Lennox.
\newblock Multivariate batch to batch optimisation of fermentation processes to
  improve productivity.
\newblock \emph{Journal of Process Control}, 108:\penalty0 148--158, 2021.

\bibitem[Blei et~al.(2017)Blei, Kucukelbir, and McAuliffe]{VI_Review}
David~M. Blei, Alp Kucukelbir, and Jon~D. McAuliffe.
\newblock Variational inference: A review for statisticians.
\newblock \emph{Journal of the American Statistical Association}, 112\penalty0
  (518):\penalty0 859--877, 2017.

\bibitem[Butcher(2000)]{ODE_Review}
{John C} Butcher.
\newblock Numerical methods for ordinary differential equations in the 20th
  century.
\newblock \emph{Journal of Computational and Applied Mathematics}, 125\penalty0
  (1--2):\penalty0 1--29, 2000.

\bibitem[Chen(2018)]{torchdiffeq}
Ricky T.~Q. Chen.
\newblock torchdiffeq, 2018.
\newblock \url{https://github.com/rtqichen/torchdiffeq}.

\bibitem[Chen et~al.(2018)Chen, Rubanova, Bettencourt, and
  Duvenaud]{chen2019_NODE}
Ricky T.~Q. Chen, Yulia Rubanova, Jesse Bettencourt, and David Duvenaud.
\newblock Neural ordinary differential equations.
\newblock \emph{Advances in Neural Information Processing Systems}, 31, 2018.

\bibitem[Cuthrell and Biegler(1989)]{cuthrell1989simultaneous}
{James E} Cuthrell and {Lorenz T} Biegler.
\newblock Simultaneous optimization and solution methods for batch reactor
  control profiles.
\newblock \emph{Computers \& Chemical Engineering}, 13\penalty0
  (1--2):\penalty0 49--62, 1989.

\bibitem[Daoutidis et~al.(2024)Daoutidis, Lee, Rangarajan, Chiang, Gopaluni,
  Schweidtmann, Harjunkoski, Mercangöz, Mesbah, Boukouvala,
  et~al.]{daoutidis2024machine}
Prodromos Daoutidis, Jay~H Lee, Srinivas Rangarajan, Leo Chiang, Bhushan
  Gopaluni, Artur~M Schweidtmann, Iiro Harjunkoski, Mehmet Mercangöz, Ali
  Mesbah, Fani Boukouvala, et~al.
\newblock Machine learning in process systems engineering: Challenges and
  opportunities.
\newblock \emph{Computers \& Chemical Engineering}, 181:\penalty0 108523, 2024.

\bibitem[Daulton et~al.(2020)Daulton, Balandat, and
  Bakshy]{daulton2020differentiable}
Samuel Daulton, Maximilian Balandat, and Eytan Bakshy.
\newblock Differentiable expected hypervolume improvement for parallel
  multi-objective bayesian optimization.
\newblock \emph{Advances in Neural Information Processing Systems},
  33:\penalty0 9851--9864, 2020.

\bibitem[Dering et~al.(2021)Dering, Swartz, and Dogan]{dering2021dynamic}
Daniela Dering, {Christopher LE} Swartz, and Neslihan Dogan.
\newblock A dynamic optimization framework for basic oxygen furnace operation.
\newblock \emph{Chemical Engineering Science}, 241:\penalty0 116653, 2021.

\bibitem[Dubois et~al.(2020)Dubois, Gordon, and Foong]{NP_Family_GitHubRepo}
Yann Dubois, Jonathan Gordon, and Andrew~YK Foong.
\newblock Neural process family, 2020.
\newblock \url{http://yanndubs.github.io/Neural-Process-Family/}.

\bibitem[Finn et~al.(2017)Finn, Abbeel, and Levine]{MAML}
Chelsea Finn, Pieter Abbeel, and Sergey Levine.
\newblock Model-agnostic meta-learning for fast adaptation of deep networks.
\newblock \emph{Proceedings of the 34th International Conference on Machine
  Learning}, 70:\penalty0 1126--1135, 2017.

\bibitem[Fogler(1999)]{fogler1999elements}
H~Scott Fogler.
\newblock \emph{Elements of chemical reaction engineering}.
\newblock Pearson Education, 1999.

\bibitem[Foong et~al.(2020)Foong, Bruinsma, Gordon, Dubois, Requeima, and
  Turner]{ConvNP}
Andrew Y.~K. Foong, Wessel~P. Bruinsma, Jonathan Gordon, Yann Dubois, James
  Requeima, and Richard~E. Turner.
\newblock Meta-learning stationary stochastic process prediction with
  convolutional neural processes.
\newblock \emph{Advances in Neural Information Processing Systems},
  33:\penalty0 8284--8295, 2020.

\bibitem[Frazier(2018)]{frazier2018tutorial}
Peter~I Frazier.
\newblock A tutorial on bayesian optimization.
\newblock \emph{arXiv 1807:02811}, 2018.

\bibitem[Garnelo et~al.(2018{\natexlab{a}})Garnelo, Rosenbaum, Maddison,
  Ramalho, Saxton, Shanahan, Teh, Rezende, and Eslami]{CNPs}
Marta Garnelo, Dan Rosenbaum, Christopher Maddison, Tiago Ramalho, David
  Saxton, Murray Shanahan, Yee~Whye Teh, Danilo Rezende, and SM~Ali Eslami.
\newblock Conditional neural processes.
\newblock \emph{International Conference on Machine Learning}, pages
  1704--1713, 2018{\natexlab{a}}.

\bibitem[Garnelo et~al.(2018{\natexlab{b}})Garnelo, Schwarz, Rosenbaum, Viola,
  Rezende, Eslami, and Teh]{Garnelo2018-NPs}
Marta Garnelo, Jonathan Schwarz, Dan Rosenbaum, Fabio Viola, Danilo~J Rezende,
  S~M~Ali Eslami, and Yee~Whye Teh.
\newblock Neural processes.
\newblock \emph{arXiv 1807:01622}, 2018{\natexlab{b}}.

\bibitem[Girolami(2008)]{bayes_inference_Differential_Eqs_Girolami2008-qi}
Mark Girolami.
\newblock Bayesian inference for differential equations.
\newblock \emph{Theoretical Computer Science}, 408\penalty0 (1):\penalty0
  4--16, 2008.

\bibitem[Huang et~al.(2020)Huang, Handel, and
  Song]{bayesian_parameter_estimation}
Hanwen Huang, Andreas Handel, and Xiao Song.
\newblock A {Bayesian} approach to estimate parameters of ordinary differential
  equation.
\newblock \emph{Computational Statistics}, 35\penalty0 (3):\penalty0
  1481--1499, 2020.

\bibitem[J\"apel and Buyel(2022)]{gp_bo_example_chromatography}
Ronald~Colin J\"apel and Johannes~Felix Buyel.
\newblock Bayesian optimization using multiple directional objective functions
  allows the rapid inverse fitting of parameters for chromatography
  simulations.
\newblock \emph{Journal of Chromatography A}, 1679:\penalty0 463408, 2022.

\bibitem[Jones et~al.(1998)Jones, Schonlau, and Welch]{Jones_EI_98}
Donald Jones, Matthias Schonlau, and William Welch.
\newblock Efficient global optimization of expensive black-box functions.
\newblock \emph{Journal of Global Optimization}, 13\penalty0 (4):\penalty0
  455--492, 1998.

\bibitem[Jorayev et~al.(2022)Jorayev, Russo, Tibbetts, Schweidtmann, Deutsch,
  Bull, and Lapkin]{Multi_Objective_BO_for_synthesis_Jorayev}
Perman Jorayev, Danilo Russo, Joshua~D. Tibbetts, Artur~M. Schweidtmann, Paul
  Deutsch, Steven~D. Bull, and Alexei~A. Lapkin.
\newblock Multi-objective {Bayesian} optimisation of a two-step synthesis of
  p-cymene from crude sulphate turpentine.
\newblock \emph{Chemical Engineering Science}, 247:\penalty0 116938, 2022.

\bibitem[Kidger(2021)]{diffrax_kidger2021on}
Patrick Kidger.
\newblock {O}n {N}eural {D}ifferential {E}quations {[PhD Thesis]. University of
  Oxford}.
\newblock 2021.

\bibitem[Kim et~al.(2019)Kim, Mnih, Schwarz, Garnelo, Eslami, Rosenbaum,
  Vinyals, and Teh]{Kim_ANP}
Hyunjik Kim, Andriy Mnih, Jonathan Schwarz, Marta Garnelo, Ali Eslami, Dan
  Rosenbaum, Oriol Vinyals, and Yee~Whye Teh.
\newblock Attentive neural processes.
\newblock \emph{arXiv 1901:05761}, 2019.

\bibitem[Kingma and Welling(2019)]{Kingma_VAE}
Diederik~P. Kingma and Max Welling.
\newblock An introduction to variational autoencoders.
\newblock \emph{Foundations and Trends in Machine Learning}, 12\penalty0
  (4):\penalty0 307--392, 2019.

\bibitem[Kocijan and Hvala(2013)]{kocijan2013sequencing}
Juš Kocijan and Nadja Hvala.
\newblock Sequencing batch-reactor control using {Gaussian}-process models.
\newblock \emph{Bioresource Technology}, 137:\penalty0 340--348, 2013.

\bibitem[Kohl et~al.(2024)Kohl, Zuo, Muir, Hornung, Polyzos, Zhu, Wang, and
  Alexander]{gp_bo_example_continous_flow}
Thomas~M. Kohl, Yan Zuo, Benjamin~W. Muir, Christian~H. Hornung, Anastasios
  Polyzos, Yutong Zhu, Xingdong Wang, and David L.~J. Alexander.
\newblock Machine-learning assisted optimisation during heterogeneous
  photocatalytic degradation utilising a static mixer under continuous flow.
\newblock \emph{Reaction Chemistry \& Engineering}, 9:\penalty0 882--893, 2024.

\bibitem[Lawrence(2003)]{lawrence2003gaussian}
Neil Lawrence.
\newblock Gaussian process latent variable models for visualisation of high
  dimensional data.
\newblock \emph{Advances in Neural Information Processing Systems}, 16, 2003.

\bibitem[Lim et~al.(1986)Lim, Tayeb, Modak, and Bonte]{lim1986computational}
HC~Lim, YJ~Tayeb, JM~Modak, and P~Bonte.
\newblock Computational algorithms for optimal feed rates for a class of
  fed-batch fermentation: Numerical results for penicillin and cell mass
  production.
\newblock \emph{Biotechnology and Bioengineering}, 28\penalty0 (9):\penalty0
  1408--1420, 1986.

\bibitem[Lyu et~al.(2018)Lyu, Xue, Yang, Yan, Hong, Zeng, and
  Zhou]{BO_for_circuit_design_Lyu}
Wenlong Lyu, Pan Xue, Fan Yang, Changhao Yan, Zhiliang Hong, Xuan Zeng, and
  Dian Zhou.
\newblock An efficient {Bayesian} optimization approach for automated
  optimization of analog circuits.
\newblock \emph{IEEE Transactions on Circuits and Systems I: Regular Papers},
  65\penalty0 (6):\penalty0 1954--1967, 2018.

\bibitem[Nguyen and Grover(2022)]{TNPs}
Tung Nguyen and Aditya Grover.
\newblock Transformer neural processes: Uncertainty-aware meta learning via
  sequence modeling.
\newblock \emph{arXiv 2207:04179}, 2022.

\bibitem[Norcliffe et~al.(2021)Norcliffe, Bodnar, Day, Moss, and
  Li{\`o}]{Norcliffe2021-NODEPs}
Alexander Norcliffe, Cristian Bodnar, Ben Day, Jacob Moss, and Pietro Li{\`o}.
\newblock Neural {ODE} processes.
\newblock \emph{arXiv 2103:12413}, 2021.

\bibitem[Patrón and Ricardez-Sandoval(2024)]{patron2024economically}
{Gabriel D} Patrón and Luis Ricardez-Sandoval.
\newblock Economically optimal operation of recirculating aquaculture systems
  under uncertainty.
\newblock \emph{Computers and Electronics in Agriculture}, 220:\penalty0
  108856, 2024.

\bibitem[Patrón et~al.(2024)Patrón, Toffolo, and
  Ricardez-Sandoval]{patron2024economic}
Gabriel~D Patrón, Kayden Toffolo, and Luis Ricardez-Sandoval.
\newblock Economic model predictive control for packed bed chemical looping
  combustion.
\newblock \emph{Chemical Engineering and Processing: Process Intensification},
  page 109731, 2024.

\bibitem[Paulson and
  Tsay(2025)]{tsay_paulson2024bayesianoptimizationflexibleefficient}
Joel~A Paulson and Calvin Tsay.
\newblock Bayesian optimization as a flexible and efficient design framework
  for sustainable process systems.
\newblock \emph{Current Opinion in Green and Sustainable Chemistry},
  51:\penalty0 100983, 2025.

\bibitem[Qing et~al.(2025)Qing, Langdon, Lee, Shafei, van~der Wilk, Tsay, and
  Misener]{SANODEP}
Jixiang Qing, Rebecca~D. Langdon, Robert~Matthew Lee, Behrang Shafei, Mark
  van~der Wilk, Calvin Tsay, and Ruth Misener.
\newblock System-aware neural {ODE} processes for few-shot {Bayesian}
  optimization.
\newblock \emph{Transactions on Machine Learning Research}, 2025.

\bibitem[Rasmussen and Williams(2006)]{GP_Rasmussen_Book}
{Carl Edward} Rasmussen and {Christopher K. I.} Williams.
\newblock \emph{Gaussian processes for machine learning}, volume~2.
\newblock MIT Press, United States, 2006.

\bibitem[Ravi and Larochelle(2017)]{Ravi2017}
Sachin Ravi and Hugo Larochelle.
\newblock Optimization as a model for few-shot learning.
\newblock \emph{International Conference on Learning Representations}, 2017.

\bibitem[San and Stephanopoulos(1989)]{profit_paper}
Ka-Yiu San and Gregory Stephanopoulos.
\newblock Optimization of fed-batch penicillin fermentation: A case of singular
  optimal control with state constraints.
\newblock \emph{Biotechnology and Bioengineering}, 34\penalty0 (1):\penalty0
  72--78, 1989.

\bibitem[Sandu et~al.(1997)Sandu, Verwer, Van~Loon, Carmichael, Potra, Dabdub,
  and Seinfeld]{Stiff_ODEs}
A.~Sandu, J.G. Verwer, M.~Van~Loon, G.R. Carmichael, F.A. Potra, D.~Dabdub, and
  J.H. Seinfeld.
\newblock Benchmarking stiff {ODE} solvers for atmospheric chemistry
  problems-{I}. implicit vs explicit.
\newblock \emph{Atmospheric Environment}, 31\penalty0 (19):\penalty0
  3151--3166, 1997.

\bibitem[Schilter et~al.(2024)Schilter, Gutierrez, Folkmann, Castrogiovanni,
  García-Durán, Zipoli, Roch, and Laino]{BO_for_auto_lab_init_conds_Schilter}
Oliver Schilter, Daniel~Pacheco Gutierrez, Linnea~M. Folkmann, Alessandro
  Castrogiovanni, Alberto García-Durán, Federico Zipoli, Loïc~M. Roch, and
  Teodoro Laino.
\newblock Combining {Bayesian} optimization and automation to simultaneously
  optimize reaction conditions and routes.
\newblock \emph{Chemical Science}, 15\penalty0 (21):\penalty0 7916--7926, 2024.

\bibitem[Shahriari et~al.(2015)Shahriari, Swersky, Wang, Adams, and
  de~Freitas]{humanoutofloop}
Bobak Shahriari, Kevin Swersky, Ziyu Wang, Ryan~P. Adams, and Nando de~Freitas.
\newblock Taking the human out of the loop: A review of {Bayesian}
  optimization.
\newblock \emph{Proceedings of the IEEE}, 104\penalty0 (1):\penalty0 148--175,
  2015.

\bibitem[Shangguan et~al.(2021)Shangguan, Lin, Wu, and Xu]{NP_BO}
Zhongkai Shangguan, Lei Lin, Wencheng Wu, and Beilei Xu.
\newblock Neural process for black-box model optimization under {Bayesian}
  framework.
\newblock \emph{arXiv 2104:02487}, 2021.

\bibitem[Shapovalova and Tsay(2025)]{SHAPOVALOVA2025469}
Mariia Shapovalova and Calvin Tsay.
\newblock Training neural {ODEs} using fully discretized simultaneous
  optimization.
\newblock \emph{IFAC-PapersOnLine}, 59\penalty0 (6):\penalty0 469--474, 2025.

\bibitem[Shin et~al.(2023)Shin, Song, Shin, Lee, and
  Seo]{gp_bo_example_catalysts}
Sangsoo Shin, Hyeongyun Song, Yeon~Su Shin, Jaegeun Lee, and Tae~Hoon Seo.
\newblock Bayesian optimization of wet-impregnated {Co-Mo/Al2O3} catalyst for
  maximizing the yield of carbon nanotube synthesis.
\newblock \emph{Nanomaterials}, 14\penalty0 (1):\penalty0 75, 2023.

\bibitem[Shokry et~al.(2018)Shokry, Vicente, Escudero, Pérez-Moya, Graells,
  and Espuña]{shokry2018data}
Ahmed Shokry, Patricia Vicente, Gerard Escudero, Montserrat Pérez-Moya,
  Moisès Graells, and Antonio Espuña.
\newblock Data-driven soft-sensors for online monitoring of batch processes
  with different initial conditions.
\newblock \emph{Computers \& Chemical Engineering}, 118:\penalty0 159--179,
  2018.

\bibitem[Slattery et~al.(2024)Slattery, Wen, Tenblad, Sanjosé-Orduna,
  Pintossi, den Hartog, and Noël]{slattery2024automated}
Aidan Slattery, Zhenghui Wen, Pauline Tenblad, Jesús Sanjosé-Orduna, Diego
  Pintossi, Tim den Hartog, and Timothy Noël.
\newblock Automated self-optimization, intensification, and scale-up of
  photocatalysis in flow.
\newblock \emph{Science}, 383\penalty0 (6681):\penalty0 eadj1817, 2024.

\bibitem[Thebelt et~al.(2022)Thebelt, Wiebe, Kronqvist, Tsay, and
  Misener]{thebelt2022maximizing}
Alexander Thebelt, Johannes Wiebe, Jan Kronqvist, Calvin Tsay, and Ruth
  Misener.
\newblock Maximizing information from chemical engineering data sets:
  Applications to machine learning.
\newblock \emph{Chemical Engineering Science}, 252:\penalty0 117469, 2022.

\bibitem[Titsias and Lawrence(2010)]{titsias2010bayesian}
Michalis Titsias and Neil~D. Lawrence.
\newblock Bayesian {Gaussian} process latent variable model.
\newblock \emph{Proceedings of the Thirteenth International Conference on
  Artificial Intelligence and Statistics}, 9:\penalty0 844--851, 2010.

\bibitem[Tsay and Baldea(2019)]{tsay2019110th}
Calvin Tsay and Michael Baldea.
\newblock 110th anniversary: {Using} data to bridge the time and length scales
  of process systems.
\newblock \emph{Industrial \& Engineering Chemistry Research}, 58\penalty0
  (36):\penalty0 16696--16708, 2019.

\bibitem[Vinyals et~al.(2016)Vinyals, Blundell, Lillicrap, and
  Wierstra]{vinyals2016matching}
Oriol Vinyals, Charles Blundell, Timothy Lillicrap, and Daan Wierstra.
\newblock Matching networks for one shot learning.
\newblock \emph{Advances in Neural Information Processing Systems}, 29, 2016.

\bibitem[Wang and Dowling(2022)]{BO_for_design_of_materials_wang}
Ke~Wang and Alexander~W. Dowling.
\newblock Bayesian optimization for chemical products and functional materials.
\newblock \emph{Current Opinion in Chemical Engineering}, 36:\penalty0 100728,
  2022.

\bibitem[Wilson et~al.(2018)Wilson, Hutter, and Deisenroth]{Wilson_BO_Acq}
James Wilson, Frank Hutter, and Marc Deisenroth.
\newblock Maximizing acquisition functions for {Bayesian} optimization.
\newblock \emph{Advances in Neural Information Processing Systems}, 31, 2018.

\bibitem[Wright and Nocedal(1999)]{wright1999numerical}
Stephen Wright and Jorge Nocedal.
\newblock \emph{Numerical optimization}.
\newblock Springer Science, 1999.

\bibitem[Xie et~al.(2025)Xie, Zhang, Paulson, and Tsay]{xie2025global}
Yilin Xie, Shiqiang Zhang, Joel Paulson, and Calvin Tsay.
\newblock Global optimization of {Gaussian} process acquisition functions using
  a piecewise-linear kernel approximation.
\newblock pages 2296--2304, 2025.

\bibitem[Xie et~al.(2026)Xie, Zhang, Qing, Misener, and Tsay]{xie2025bogrape}
Yilin Xie, Shiqiang Zhang, Jixiang Qing, Ruth Misener, and Calvin Tsay.
\newblock {BoGrape: Bayesian} optimization over graphs with shortest-path
  encoded.
\newblock \emph{arXiv 2503:05642}, 2026.

\bibitem[Zhou et~al.(2009)Zhou, Li, Qian, Chen, and
  Kraslawski]{zhou2009optimizing}
Hua Zhou, Xiuxi Li, Yu~Qian, Yun Chen, and Andrzej Kraslawski.
\newblock Optimizing the initial conditions to improve the dynamic flexibility
  of batch processes.
\newblock \emph{Industrial \& Engineering Chemistry Research}, 48\penalty0
  (13):\penalty0 6321--6326, 2009.

\bibitem[Zhou et~al.(2015)Zhou, Chen, and Song]{zhou2015recursive}
Le~Zhou, Junghui Chen, and Zhihuan Song.
\newblock Recursive {Gaussian} process regression model for adaptive quality
  monitoring in batch processes.
\newblock \emph{Mathematical Problems in Engineering}, 2015:\penalty0 761280,
  2015.

\end{thebibliography}
\bibliographystyle{plainnat} 


\end{document}